 \documentclass[final,3p,times]{elsarticle}
\usepackage{amssymb}
\usepackage{amsmath}
\usepackage{amsthm}
\usepackage[caption=false]{subfig}
\usepackage{comment}
\usepackage{graphicx}
\usepackage{booktabs}
\usepackage[breaklinks]{hyperref}
\theoremstyle{plain}
\newtheorem{theorem}{Theorem}[section]
\newtheorem{lemma}[theorem]{Lemma}

\theoremstyle{definition}
\newtheorem{definition}[theorem]{Definition}
\newtheorem{example}[theorem]{Example}
\newtheorem*{solution}{Solution}

\theoremstyle{remark}
\newtheorem{remark}{Remark}

\journal{Optical and Quantum Electronics}

\begin{document}

\begin{frontmatter}



\title{On a Class of Multi-Dimensional Non-linear Time-Fractional Fokker-Planck Equations Capturing Brownian Motion}


\author[1,*]{Neetu Garg}
\author[1]{Varsha R}
\affiliation[1]{organization={Department of Mathematics, National Institute of Technology Calicut},
            city={Kozhikode},
            postcode={673601}, 
            state={Kerala},
            country={India}}
\affiliation[*]{Corresponding author Email: neetu@nitc.ac.in}
\begin{abstract}
The time-fractional Fokker-Planck equation is a key model for characterizing anomalous diffusion, stochastic transport, and non-equilibrium statistical mechanics with applications in finance, chaotic dynamics, optical physics, and biological systems. In this work, we develop a semi-analytical solution for the multi-dimensional time-fractional Fokker-Planck equation employing the Laplace residual power series method. This method blends the Laplace transform and the traditional residual power series method, guaranteeing efficient solutions incorporating the memory and nonlocal effects. To validate the accuracy and effectiveness of the approach, we address several examples including non-linear problems in multi-dimensions, and analyze the evolution of errors. The numerical simulations are compared with existing methods to confirm the adopted method's strength. The smooth and stable error evolution promises that the suggested method is a powerful tool for analyzing time-fractional Fokker-Planck equations. 
\end{abstract}

\begin{highlights}
\item The Fokker-Planck equation is a significant partial differential equation describing the evolution of the probability density function and fractional Brownian motion. 
\item This study deals with the Caputo time fractional Fokker-Planck equation incorporating memory effect and non-local property, which depicts anomalous diffusion, stochastic transport,
and non-equilibrium statistical mechanics with applications in finance, chaotic dynamics, optical physics, etc.
\item  The hybrid technique Laplace residual power series method, is employed to construct the solution. 
\item An approximate series solution method is formulated for numerical examples, encompassing both linear and nonlinear multidimensional problems.
\item  The numerical simulations match up with established methods that validate the efficacy of the proposed approach.
\end{highlights}


\end{frontmatter}


\section{Introduction}
\label{sec1}

Fractional calculus is an extension of classical calculus concerning the study of integrals and derivatives of arbitrary order. During the last three decades, fractional calculus has become popular and interesting among researchers. Due to the wide implications and applications of fractional operators, many complex real-life phenomena are described using them. Numerous complex systems seen in the fields of physics, biology, and social-economics exhibit dynamics that cannot be sufficiently characterized by traditional integer-order models. Thus, they have been altered into fractional models for better design of complex phenomena. This is because of the privilege of choosing a suitable fractional order for the model and the memory effect property held by the fractional operators. The fractional order derivatives are non-local, because they characterize the long-term memory of the system \cite{KBO74, IP99}. 
 The application of fractional calculus includes modeling human diseases, signal processing, fluid dynamics, control theory, population dynamics, etc. \cite{H00, M11, RTVV11}.

Generally, differential equations having fractional order derivatives are called the fractional differential equations (FDE). In this study, we aim to analyze one of the most important non-linear FDEs called time-fractional Fokker-Planck equations (TFFPEs). Adriaan Fokker and Max Planck developed the classical Fokker-Planck equation (FPE) to characterize the Brownian motion of particles and the time evolution of the probability density function of a particle’s position \cite{R12}. The FPE governs the 
statistical randomness that occurs in the Brownian motion of particles due to random forces. The FPE is a widely used equation in statistical physics that characterizes the complicated dynamics observed in complex physical, biological, and physio chemical systems. In addition, it emerges in modeling various natural phenomena such as quantum optics, probability flux, theoretical biology, circuit theory, etc. \cite{DHM20, EG16}. The general form of one-dimensional classical FPE representing the motion of particle's concentration  in the spatial variable $\zeta$ at time $\tau$ has the following form:
\begin{equation}
    \frac{\partial }{\partial t} \nu(\zeta, \tau)= \left(- \frac{\partial}{\partial \zeta} f_1(\zeta )+   \frac{\partial ^2 }{\partial \zeta ^2} f_2(\zeta)\right)  \nu(\zeta, \tau), 
\end{equation}
with the initial condition $v(\zeta,0)= \varphi (\zeta), ~\zeta \in \mathbb{R}.$ 
Here, $f_1(\zeta)$ and $f_2(\zeta)$ are non-negative smooth functions called drift and diffusion coefficients, respectively. The scenario in which the drift and diffusion coefficients are functions of time $\tau$ and $\nu$ will yield the nonlinear Fokker-Planck equation. These nonlinear Fokker-Planck equations characterize numerous scientific phenomena, including polymer physics, pattern generation, surface physics, etc.

Recently, the integer order time derivative in the Fokker-Planck equation has been substituted with the fractional order derivative. Subsequently, numerous mathematicians began to exhibit significant interest in the research of TFFPEs. This equation arose from the necessity to simulate anomalous diffusion, when conventional diffusion diverges from Brownian motion. The ability to characterize anomalous diffusion in continuous time random walks, where particles adhere to a random time interval between jumps, renders the TFFPEs more significant than the classical FPE. The TFFPE which describes the motion of concentration of particles $\nu(\zeta,\ \tau)$ in one-dimensional space and at time $\tau$ has the general form
\cite{HAS24},
\begin{equation} \label{eq:1.4}
	D_{\tau}^{\gamma} \nu(\zeta, \tau)=-\frac{\partial}{\partial \zeta} \left(f_{1}(\zeta, \tau,\nu)\nu(\zeta, \tau)\right)+\frac{\partial^{2}}{\partial \zeta^{2}}\left( f_{2}(\zeta, \tau,\nu)\nu(\zeta, \tau)\right) , \ \zeta, \ \tau \geq 0, \end{equation}
with initial condition
\begin{equation} \label{eq:1.6} 
	\nu(\zeta,0)=\varphi(\zeta),  \ \zeta\in \mathbb{R}.
    \end{equation}
     The $D_{\tau}^{\gamma}$ indicates the Caputo time-fractional derivative with order $\gamma$ ($0<\gamma\leq 1).$ For $\gamma=1,$ the TFFPE simplifies to the traditional FPE. The fractional derivative in the TFFPE functions as a potent instrument for connecting microscopic stochastic processes to emerging macroscopic statistical phenomena, a fundamental concept in statistical mechanics. The TFFPEs have numerous applications in physical problems such as anomalous diffusion in physio-chemical systems, polymer dynamics, intracellular transport in biological systems, and so on \cite{WWZMC23}. Although it has a wider range of applications, the methods for obtaining the exact solution are very limited. Recently, K. Singla and N. Leonenko \cite{SL25} utilized the Lie group classification to reduce the fractional FPE into a fractional ordinary differential equation and thereby obtained the exact solution as a power series. In addition, many analytical and numerical approaches have been proposed to extract approximate solutions, including neural network method, residual power series method, fractional reduced differential transform method, the fractional homotopy perturbation transform method, finite difference method, mixed finite element method, and so on \cite{WWLZ22, YKK15,  DC15, K13, KMA24, WCK25}.
     
     The main concern of this study is to analyze and investigate the approximate semi-analytical solution of the nonlinear TFFPE in multidimensional space having the general form:
      \begin{equation} \label{EQ4}
 	D_\tau^ \gamma \nu(\vec{\zeta},\ 
    \tau)=\left[-\sum_{i=1}^{d}\frac{\partial}{\partial \zeta _i} f_{1,i}(\vec{\zeta}, \ \tau,\nu)+\sum_{i=1}^{d}\sum_{j=1}^{d}\frac{\partial ^2}{\partial \zeta_i \zeta_j}f_{2,i,j}(\vec{\zeta}, \ \tau, \ \nu)\right]\nu(\vec{\zeta}, \ \tau),
 \end{equation} where $\vec{\zeta} = (\zeta_1, \zeta_2,\hdots ,\zeta_d)\in\mathbb{R}^d,$ $d=1,2,3.$ 
  The literature on multidimensional TFFPE is limited and somewhat constrained, especially compared to the integer order condition. This arises from the mathematical intricacies and numerical difficulties. Thus, obtaining a solution for multidimensional TFFPE is an attractive problem among researchers. The motive of this manuscript is to explore the application of the Laplace residual power series (LRPS) method to solve the multidimensional TFFPEs. The LRPS is a semi-analytical method formulated by the blend of the Laplace transform (LT) and the traditional residual power series method (RPSM) \cite{AEM15}. In 2020, Eriqat et al. \cite{EEOAM20} introduced this novel hybrid technique by promoting RPSM using the Laplace transform (LT). The application of LT reduces the FDE to a simple algebraic equation. Then, we exploit the residual function and some limit conditions to compute the coefficients in the power series. Finally, applying the inverse Laplace operator provides an approximate series solution for the equation considered. 
  
  The LRPS approach was first investigated to address non-linear pantograph equations \cite{EEOAM20}. Recently, employing this technique, El-Ajou et al. computed the series solution for one-dimensional time-fractional Schrödinger equation \cite{ESDQA24}, Pant et al. solved the 2-dimensional FDEs \cite{PASE24}, and Oqielat et al. established the series solution for the time-fractional reaction-diffusion model \cite{OEEAA23}. Few fractional models analyzed using this technique are temporal fractional  Newell–Whitehead–Segel equations \cite{LN23},  fractional reaction-diffusion Brusselator model \cite{ASIAN22}, time-fractional coupled Boussinesq-Burger equations \cite{SBSA22}, time-fractional Black-Scholes equations \cite{LO23}, Bacteria growth model \cite{OEAOEH23}, SIR epidemic model \cite{QS23}, 
   and so on \cite{AEAES23, E21}.
   
   The biggest advantage of the proposed approach is the absence of the need for linearization or perturbation to address the nonlinearities within the systems. Furthermore, the solution is devoid of round-off error, often seen in numerical methods, since discretization is not used. Therefore, it is intriguing to use this strategy due to its simplicity, effectiveness, and dependability. In this article, we address multi-dimensional TFFPE \eqref{EQ4} using the LRPS technique to provide a precise approximate solution. We examine and evaluate the acquired solutions and their corresponding errors using tables and graphs. To assess the efficacy of the suggested strategy, we juxtapose our results with the current literature. Moreover, we analyze the impact of a control term in the numerical simulation of TFFPE.

The paper is framed as follows: The second section reviews the properties of the Caputo fractional derivative and LT. In the third section, the outline of the LRPS method is explained. A few examples are solved using the LRPSM in the fourth section. To convey the capability and precision of the adopted method, numerical results are depicted in the fifth section. In the sixth section, we discuss the consequence of the control term in the numerical simulation of TFFPE. We conclude in the last section.

\section{Preliminaries}\label{sec2}

This section describes some basics of fractional calculus, the Laplace transform, and some related required theorems.

\begin{definition}
\cite{IP99}	The Caputo fractional derivative of $f(\tau)$ of order $\gamma \in \mathbb{R}^{+}$ is given as:
	
	\begin{equation}
		D_{\tau}^{\gamma} f(\tau)=\left\{\begin{array}{lr}
			\frac{1}{\Gamma(l-\gamma)} \int_{0}^{\tau} \frac{f^{(l)}(\rho)}{(\tau-\rho)^{\gamma+1-l}} d \rho, & l-1 <\gamma<l \\
			f^{(l)}(\tau), & \gamma=l \in \mathbb{N}
		\end{array}\right.
	\end{equation}
	where $\Gamma(.)$ is the Gamma function.
\end{definition}

\begin{definition} \cite{IP99}
	The Mittag-Leffler function is defined as:
	\begin{equation}
		E_\gamma(\chi)=\sum_{n=0}^\infty \frac{\chi^n}{\Gamma(n\gamma+1)}, \hspace{0.2cm} \gamma\in \mathbb{C}, \hspace{0.1cm} Re(\gamma)>0.
	\end{equation}
\end{definition}

\begin{lemma}  \label{lem1}
       \cite{IP99} Consider a piecewise continuous function $\nu(\zeta, \tau)$ defined in the domain $I \times[0 \times \infty)$ of the exponential order $\delta$. Assume $\mathcal{L}\{\nu(\zeta, \tau)\}=\mathcal{V}(\zeta, s)$, then

       \begin{enumerate}
       	\item $ \mathcal{L}\left\{D_{\tau}^{\gamma} \nu(\zeta, \tau)\right\}=s^{\gamma} \mathcal{V}(\zeta, s)-\sum_{i=0}^{l-1} s^{\gamma-i-1} D_{\tau}^{i} \nu(\zeta, 0), \  \text{for } l-1<\gamma \leq l,$ $ l \in \mathbb{N} .$
       	\item $ \mathcal{L}\left\{D_{\tau}^{n \gamma} \nu(\zeta, \tau)\right\}=s^{n \gamma} \mathcal{V}(\zeta, s)-\sum_{j=0}^{n-1} s^{(n-j) \gamma-1} D_{\tau}^{j \gamma} \nu(\zeta, 0),$ \ $ 0<\gamma<1 $, where $D_{\tau}^{n \gamma}=D_{\tau}^{\gamma} \cdot D_{\tau}^{\gamma} \ldots D_{\tau}^{\gamma}$  ($n$  times), $n\in \mathbb{N}$. 
       	\item $\lim _{s \rightarrow \infty} s \mathcal{V}(\zeta, s)=\nu(\zeta, 0) .$
       \end{enumerate}
\end{lemma}

\begin{definition}
	A fractional power series (FPS) about $\tau=\tau_{0}$ is defined as follows:
	
	$$
	\begin{aligned}
		\sum_{n=0}^{\infty} p_{n}(\zeta)\left(\tau-\tau_{0}\right)^{n \gamma},
		l-1 < \gamma<l, l & \in \mathbb{N}, \tau \geq \tau_{0}.
	\end{aligned}
	$$
\end{definition}

\begin{theorem}
	\cite{E21} Let $\mathcal{V}(\zeta,s)=\mathcal{L}\{\nu(\zeta,\tau)\}$ and suppose that $\mathcal{V}(\zeta,s)$ has the $\gamma$ -singular Laurent series representation as $\mathcal{V}(\zeta,s)=\sum_{n=0}^\infty \frac{p_n(\zeta)}{s^{n\gamma+1}}$ where $0<\gamma\leq 1$. 
	Then, $$
	p_n(\zeta)=D_\tau^{n\gamma}\nu(\zeta,0).$$ 
\end{theorem}

\begin{remark}
	Suppose $\mathcal{V}(\zeta,s)=\sum_{n=0}^\infty \frac{p_n(\zeta)}{s^{n\gamma+1}},$ then inverse LT of $\mathcal{V}(\zeta,s)$ is given as, 
	\begin{equation*}
		\nu(\zeta,\tau)=\sum_{n=0}^\infty \frac{D_\tau^{n\gamma}\nu(\zeta,0)}{\Gamma(n\gamma+1)}\tau^{n\gamma}.
	\end{equation*}
	This is the FPS representation of $\nu(\zeta,\tau)$ about $\tau=0$.
	
\end{remark}

The main concept related to the series solution is the convergence of the series. A necessary theorem for the convergence of FPS is stated below: 
\begin{theorem} \cite{E21}
    The  FPS $\sum_{n=0}^\infty p_n(\zeta)(\tau-\tau_0)^{n\gamma}$ about $\tau=\tau_0$ have the following three possible ways of convergence:
    \begin{enumerate}
        \item[(i)] The series converges only at $\tau=\tau_0,$ \textit{i.e.} radius of convergence is zero.
        \item[(ii)] The series converges for all $\tau \geq \tau_0$, \textit{i.e.} the series have an infinite radius of convergence.
        \item[(iii)] There exists a real number $R>0,$ such that the series converges within the range $\tau_0 \leq \tau \leq \tau_0 +R$ and diverges otherwise. The quantity $R$ is the radius of convergence.
    \end{enumerate}
\end{theorem}
\vskip 0.5cm
\begin{theorem}
	\cite{E21}
	Let $\nu(\zeta,\tau)$ is a piecewise continuous function on $I\times [0,\infty)$ having exponential order $\delta$ and $\mathcal{V}(\zeta,s)=\mathcal{L}\{\nu(\zeta,\tau)\}.$ Suppose $\mathcal{V}(\zeta,s)$ has the $\gamma$-singular Laurent series expression as $\mathcal{V}(\zeta,s)= \sum_{n=0}^{\infty}\frac{p_n(\zeta)}{s^{n\gamma+1}}$ and it satisfy $\left |s\mathcal{L}\left\{D_\tau^{(n+1)\gamma}\nu(\zeta,\tau)\right\} \right|\leq \mathcal{M}(\zeta)$ for $\gamma \in (0,1] $ on $I \times [0,\infty)$. Then the remainder $R_n(\zeta,s)$ of the $\gamma$-singular Laurent series of $\mathcal{V}(\zeta,s)$ have a bound such that,
	$\left|R_n(\zeta,s) \right|\leq \frac{\mathcal{M}(\zeta)}{s^{(n+1)\gamma+1}}$ for $\zeta \in I$.
\end{theorem}

\section{Construction of LRPS solution for TFFPE} 

This section discusses the methodology to construct an LRPS solution for TFFPE. The primary concept of this technique is to transform the FDE into a simple algebraic equation and use certain limit conditions to evaluate the coefficients in the approximate series solution.

Transforming Eq. (\ref{eq:1.4}) by applying the Laplace operator and invoking the initial condition, we get
\begin{equation} \label{eq:3.1}
	\mathcal{V}(\zeta, s)=\frac{\varphi(\zeta)}{s}+\frac{1}{s^{\gamma}} \mathcal{L}\left\{\left[-\frac{\partial}{\partial \zeta} f_{1}(\zeta, \tau,\nu)+\frac{\partial^{2}}{\partial \zeta^{2}} f_{2}(\zeta, \tau,\nu)\right] \mathcal{L}^{-1}\{\mathcal{V}(\zeta, s)\}\right\},
\end{equation} where $\mathcal{V}(\zeta,s)=\mathcal{L}\{\nu(\zeta,\tau)\}.$

The LRPS method proposes the solution of TFFPE in the series form as follows, 
\begin{equation}
	\mathcal{V}(\zeta, s)=\sum_{n=0}^{\infty} \frac{p_{n}(\zeta)}{s^{n \gamma+1}},\hspace{0.2cm} \zeta, \tau \geq 0, \hspace{0.2cm} s>0.
\end{equation}
Now, we express the $k^{t h}$ truncated series of $\mathcal{V}(\zeta, s)$ as,

\begin{equation}
	\mathcal{V}_{k}(\zeta, s)=\sum_{n=0}^{k} \frac{p_{n}(\zeta)}{s^{n \gamma+1}}, \hspace{0.2cm} \zeta, \tau \geq 0,\hspace{0.2cm} s>0.
\end{equation}

Next, the Laplace residual function (LRF) is defined as,
\begin{equation}
	\mathcal{L} \mathcal{R}(\zeta, s)=\mathcal{V}(\zeta, s)-\frac{\varphi(\zeta)}{s}-\frac{1}{s^{\gamma}} \mathcal{L}\left\{\left[-\frac{\partial}{\partial \zeta} f_{1}(\zeta, \tau,\nu)+\frac{\partial^{2}}{\partial \zeta^{2}} f_{2}(\zeta, \tau,\nu)\right] \mathcal{L}^{-1}\{\mathcal{V}(\zeta, s)\}\right\},
\end{equation}

and the $k^{\text {th }}$ truncated LRF by replacing $\mathcal{V}(\zeta,s)$ by $\mathcal{V}_k(\zeta,s)$ in LRF as,
\begin{equation}
	\mathcal{LR}_k(\zeta, s)=\mathcal{V}_{k}(\zeta, s)-\frac{\varphi(\zeta)}{s}-\frac{1}{s^{\gamma}} \mathcal{L}\left\{\left[-\frac{\partial}{\partial \zeta} f_{1}(\zeta, \tau,\nu)+\frac{\partial^{2}}{\partial \zeta^{2}} f_{2}(\zeta, \tau,\nu)\right] \mathcal{L}^{-1}\{\mathcal{V}_{k}(\zeta, s)\}\right\}.
\end{equation}
\noindent The significant features of LRF for the construction of the LRPS solution are listed below:

\textbf{P1}. $\mathcal{LR} (\zeta, s)=0, \hspace{0.2cm} \zeta, \tau \geq 0,\hspace{0.2cm} s>0$.

\textbf{P2}. $\lim _{k \rightarrow \infty} \mathcal{L} \mathcal{R} _{k}(\zeta, s)=\mathcal{L} \mathcal{R} (\zeta, s), \hspace{0.2cm} \zeta, \tau \geq 0, \hspace{0.2cm} s>0$.

\textbf{P3}. If $\lim _{s \rightarrow \infty} s\mathcal{LR}  (\zeta, s)=0,$ then $\lim _{s \rightarrow \infty} s \mathcal{L} \mathcal{R} _{k}(\zeta, s)=0$.

\textbf{P4}. $\lim _{s \rightarrow \infty} s^{k \gamma+1} \mathcal{L} \mathcal{R} (\zeta, s)=\lim _{s \rightarrow \infty} s^{k \gamma+1} \mathcal{L} \mathcal{R} _{k}(\zeta, s)=0$, for $k=1,2,3, \ldots$ 

Using the property P4, we determine the unknown coefficient $p_n(\zeta)$. By substituting these coefficients, we derive the approximate solution $\mathcal{V}_k(\zeta,s)$ for Eq. (\ref{eq:3.1}).
The final step is to apply inverse LT to $\mathcal{V}_{k}(\zeta, s)$ to obtain an approximate series solution $\nu_{k}(\zeta,\tau)$ for the concerned TFFPE. 

\begin{remark}

For higher dimensions, such as $d = 2 , 3$, the series solution assumes the following structure:
\begin{equation}
	\mathcal{V}(\vec{\zeta},s)=\sum_{n=0}^\infty \frac{p_n(\vec{\zeta})}{s^{n\gamma+1}}, \hspace{0.2cm} \vec{\zeta} \in \mathbb{R}^d, \hspace{0.2cm} d=2,3.
\end{equation} The remaining procedure for deriving the solution follows the same procedure as in the one-dimensional case.
    
\end{remark}

\section{Applications}

In this section, we investigate a few numerical problems in multi-dimensional space to present the accuracy of the LRPS method.
\begin{example}  \label{ex1}
	Consider the linear TFFPE in one-dimension with $f_1(\zeta,\tau,\nu)=-1$ and $f_2(\zeta,\tau,\nu)=1.$
	\begin{equation} \label{eq:29}
		D_{\tau}^{\gamma} \nu(\zeta, \tau)=D_\zeta \nu(\zeta, \tau)+D_{\zeta\zeta}\nu(\zeta, \tau), \hspace{0.2cm} \zeta\in \mathbb{R},\hspace{0.2cm} \tau > 0, \hspace{0.2cm}\gamma \in(0,1],   
	\end{equation}  
    \begin{equation} \label{eq:30}
		\nu(\zeta,0)=\zeta,   
	\end{equation}  where $D_{\zeta}=\frac{\partial}{\partial \zeta}$ and $D_{\zeta \zeta}=\frac{\partial^{2}}{\partial \zeta^{2}}$ respectively. The exact solution is $\zeta+\frac{\tau^\gamma}{\Gamma(\gamma+1)}.$
    \begin{solution}
        
   By taking the LT on 
   Eq. (\ref{eq:29}) and simplifying 
   \begin{equation}  \label{eq:32}
		\mathcal{V}(\zeta,s)=\frac{\zeta}{s}+\frac{1}{s^\gamma}D_\zeta\mathcal{V}(\zeta,s)+\frac{1}{s^\gamma}D_{\zeta\zeta}\mathcal{V}(\zeta,s).\end{equation} Now, we define the LRF as \begin{equation} \label{eq:33}
		\mathcal{LR}(\zeta,s)=\mathcal{V}(\zeta,s)-\frac{\zeta}{s}-\frac{1}{s^\gamma}D_\zeta \mathcal{V}(\zeta,s)-\frac{1}{s^\gamma}D_{\zeta\zeta}\mathcal{V}(\zeta,s),\end{equation}
	and the $k^{th}$  LRF as, \begin{equation} \label{eq:34}
		\mathcal{LR}_k(\zeta,s)=\mathcal{V}_k(\zeta,s)-\frac{\zeta}{s}-\frac{1}{s^\gamma}D_\zeta \mathcal{V}_k(\zeta,s)-\frac{1}{s^\gamma}D_{\zeta\zeta}\mathcal{V}_k(\zeta,s).\end{equation}The fractional power series solution of Eq. (\ref{eq:32}) is \begin{equation}
		\mathcal{V}(\zeta,s)=\sum_{n=0}^{\infty} \frac{p_{n}(\zeta)}{s^{n \gamma+1}},
	\end{equation} where $p_0(\zeta)=\zeta$ by  Lemma \ref{lem1}. 
	Now, the $k^{th}$ truncated series is given as, 
	\begin{equation} \label{eq:35} 
		\mathcal{V}_k(\zeta,s)=\frac{\zeta}{s}+\sum_{n=1}^{k} \frac{p_{n}(\zeta)}{s^{n \gamma+1}}.
	\end{equation}
    Substituting Eq. (\ref{eq:35}) in Eq. (\ref{eq:34}) and using property \textbf{P4}, we get the unknown coefficients values as, $$p_1(\zeta)=1 \hspace{0.2cm}\text{and } p_i(\zeta)=0\hspace{0.2cm}\text{for}\hspace{0.2cm} i=2,3,...$$
	Thus \begin{equation}    
		\mathcal{V}(\zeta,s)=\frac{\zeta}{s}+\frac{1}{s^\gamma+1}.\end{equation}
	Finally by applying inverse LT, we achieve $\nu(\zeta,\tau)=\zeta+\frac{\tau^\gamma}{\Gamma(\gamma+1)}.$ In case of $\gamma=1$,  $\nu(\zeta,\tau)$ reduces to the exact solution $\zeta+\tau$.
     \end{solution}
\end{example}

\begin{example}   \label{ex2}
	Consider the linear one-dimensional TFFPE, 
	\begin{equation} \label{eq:38}
		D_\tau^{\gamma}\nu(\zeta,\tau)=-D_\zeta\bigl(\zeta\nu(\zeta,\tau)\bigr)+D_{\zeta\zeta}\left(\frac{\zeta^2}{2}\nu(\zeta,\tau)\right),\hspace{0.2cm} \zeta\in\mathbb{R},\hspace{0.2cm} \tau \geq 0,\hspace{0.2cm}\gamma \in(0,1],\end{equation}
	 \begin{equation} \label{eq:39} \nu(\zeta,0)=\zeta.\end{equation} 
	 The exact solution is $\nu(\zeta,\tau)=\zeta E_\gamma(\tau^\gamma).$ 
    \begin{solution}
        Taking the LT on Eq. (\ref{eq:38}), we have \begin{equation} \label{eq:40}
		\mathcal{V}(\zeta,s)=\frac{\zeta}{s}-\frac{1}{s^\gamma}D_\zeta(\zeta\mathcal{V}(\zeta,s))+\frac{1}{s^\gamma}D_{\zeta\zeta}\left(\frac{\zeta^2}{2}\mathcal{V}(\zeta,s)\right).\end{equation}
	Adhering to the same procedures as demonstrated in the preceding example, the LRF and $k^{th}$ LRF are in the following form, 
	\begin{equation}
		\mathcal{LR} (\zeta, s)=\mathcal{V}(\zeta, s)-\frac{\zeta}{s}+\frac{1}{s^\gamma}D_\zeta(\zeta\mathcal{V}(\zeta,s))-\frac{1}{s^\gamma}D_{\zeta\zeta}\left(\frac{\zeta^2}{2}\mathcal{V}(\zeta,s)\right),
	\end{equation}
	\begin{equation} \label{eq:44}
		\mathcal{LR}  _k(\zeta, s)=\mathcal{V}_k(\zeta, s)-\frac{\zeta}{s}+\frac{1}{s^\gamma}D_\zeta\left(\zeta\mathcal{V}_k(\zeta,s)\right)-\frac{1}{s^\gamma}D_{\zeta\zeta}\left(\frac{\zeta^2}{2}\mathcal{V}_k(\zeta,s)\right).
	\end{equation}
	
\noindent	Applying property \textbf{P4}, we yields the value of unknown coefficients as, 
	$$p_i(\zeta)=\zeta,\hspace{0.2cm} 1\leq i\leq 8.$$
	Hence, \begin{equation} \label{eq:45}
		\mathcal{V}_8(\zeta, s)=\sum_{k=0}^{8}\frac{\zeta}{s^{k\gamma +1}}.
        \end{equation}
        The final step is to apply the inverse LT on Eq. (\ref{eq:45}) to get the solution as,\begin{equation} \label{eq:47}
\nu_8(\zeta, \tau)=\zeta \sum_{k=0}^{8}\frac{\tau^ {k\gamma}}{\Gamma(k\gamma+1)}. 
\end{equation}
	In general, \begin{equation} \label{eq:46}
		\mathcal{V}(\zeta,s)=\ \sum_{k=0}^{\infty}\frac{\zeta}{s^{k\gamma +1}}. 
        \end{equation}
        Utilizing the inverse LT on Eq. (\ref{eq:46}) gives, \begin{equation} \label{eq:48}
		\nu(\zeta,\tau)=
        \sum_{k=0}^{\infty} \frac{ \zeta\ \tau^ {k\gamma}}{\Gamma(k\gamma+1)}. 
        \end{equation}
	Thus the closed-form solution of Eq. (\ref{eq:38})-(\ref{eq:39}) is $\nu(\zeta,\tau)=\zeta E_\gamma(\tau^\gamma).$ In case of $\gamma=1,$ it reduces to $\nu(\zeta,\tau)=\zeta \exp (\tau).$
    
    \end{solution}
\end{example}

\begin{example}   \label{ex4}
 Consider the non-linear one-dimensional TFFPE 
 \begin{equation}  \label{eq:21}
 	D_{\tau}^{\gamma} \nu(\zeta, \tau)=-D_\zeta\left( 3\nu^2(\zeta, \tau)-\frac{\zeta}{2}\nu(\zeta, \tau)\right)+D_{\zeta\zeta}\left(\zeta\nu^2(\zeta, \tau)\right), \hspace{0.2cm} \zeta\in\mathbb{R},\hspace{0.2cm} \tau \geq 0, \hspace{0.2cm}\gamma \in(0,1],
 \end{equation}
\begin{equation} \label{eq:22}
 	\nu(\zeta,0) =\zeta.  
 \end{equation} The exact solution of  Eqs. (\ref{eq:21})-(\ref{eq:22}) is $\nu(\zeta,\tau)=\zeta E_\gamma(\tau^\gamma).$\begin{solution}
  To obtain the solution, first we apply the LT to Eq. (\ref{eq:21}), 
  \begin{equation}   \label{eq:23}
 	\mathcal{V}(\zeta,s)=\frac{\zeta}{s}-\frac{1}{s^\gamma}D_\zeta\Bigl({3\mathcal{L}\{(\mathcal{L}^{-1}\{\mathcal{V}(\zeta,s)\})^2\}}-\frac{\zeta}{2}\mathcal{V}(\zeta,s)\Bigr)+\frac{1}{s^\gamma}D_{\zeta\zeta}\Bigl(\zeta\mathcal{L}\{(\mathcal{L}^{-1}\{\mathcal{V}(\zeta,s)\})^2\}\Bigr).\end{equation}
 \\
Using the limit condition for $\mathcal{LR}_k(\zeta,s)$ recursively, we obtain $p_i(\zeta)=\zeta$ for $1\leq i\leq 8.$  
  Thus, \begin{equation} \label{eq:26}
 	\mathcal{V}_8(\zeta,s)=
    \sum_{k=0}^8 \frac{\zeta}{s^{k\gamma+1}}. 
    \end{equation} Hence  applying inverse LT to Eq. (\ref{eq:26}), we get \begin{equation} \label{eq:27}
 	\nu_8(\zeta,\tau)=
    \sum_{k=0}^8 \frac{\zeta \tau^ {k\gamma}}{\Gamma(k\gamma+1)}.
\end{equation} So in general, \begin{equation} \label{eq:28}
 	\nu(\zeta,\tau)=
    \sum_{k=0}^\infty \frac{\zeta \tau^ {k\gamma}}{\Gamma(k\gamma+1)}.
    \end{equation}
 The closed-form solution of Eq. (\ref{eq:21})-(\ref{eq:22})  is $\zeta E_\gamma(\tau^\gamma)$
which converges to the exact solution $\zeta \exp ({\tau})$ for $\gamma=1$.
     
 \end{solution}
\end{example}
\begin{example}   \label{ex5}
Consider the linear TFFPE in 2-dimension,
\begin{equation} \label{51}
	\begin{split}
		D_{\tau}^{\gamma} \nu(\zeta_1, \zeta_2, \tau)
		=&-D_{\zeta_1}
		(\zeta_1\nu(\zeta_1,\zeta_2,\tau))
		-D_{\zeta_2}(5\zeta_2\nu(\zeta_1,\zeta_2,\tau))
		+D_{\zeta_1\zeta_1}(\zeta_1^2\nu(\zeta_1,\zeta_2,\tau))
		+D_{\zeta_1\zeta_2}\nu(\zeta_1,\zeta_2,\tau)\\
		&
		+D_{\zeta_2\zeta_1}\nu(\zeta_1,\zeta_2,\tau)+D_{\zeta_2\zeta_2}(\zeta_2^2\nu(\zeta_1,\zeta_2,\tau)),
	\end{split}
\end{equation}	
\begin{equation}
\nu(\zeta_1,\zeta_2,0)=\zeta_1, \hspace{0.2cm} (\zeta_1, \zeta_2) \in\mathbb{R}^2,\hspace{0.2cm} \tau \geq 0. \label{52}
\end{equation} The exact solution is $\nu(\zeta_1,\zeta_2,\tau)=\zeta_1E_\gamma(\tau^\gamma).$\begin{solution}
   Operating the LT on Eq. (\ref{51}) and utilising (\ref{52}), we have \begin{equation} \label{53}
	\begin{split}
	\mathcal{V}(\zeta_1,\zeta_2,s)&=\frac{\zeta_1}{s}-\frac{1}{s^\gamma}\Bigg[D_{\zeta_1}(\zeta_1\mathcal{V}(\zeta_1,\zeta_2,s))+D_{\zeta_2}(5\zeta_2\mathcal{V}(\zeta_1,\zeta_2,s))-D_{\zeta_1\zeta_1}(\zeta_1^2\mathcal{V}(\zeta_1,\zeta_2,s))
	\\&-D_{\zeta_1\zeta_2}\mathcal{V}(\zeta_1,\zeta_2,s)-D_{\zeta_2\zeta_1}\mathcal{V}(\zeta_1,\zeta_2,s)-D_{\zeta_2\zeta_2}(\zeta_2^2\mathcal{V}(\zeta_1,\zeta_2,s))\bigg].
\end{split}
\end{equation}  Here, $\mathcal{V}(\zeta_1,\zeta_2,s)$ has the following representation:
\begin{equation}
	\mathcal{V}(\zeta_1,\zeta_2,s)=\sum_{n=0}^{\infty}\frac{p_n(\zeta_1,\zeta_2)}{s^{n\gamma+1}}.
\end{equation}
Now, we employ property \textbf{P4} to evaluate the unknown coefficients of the series solution. Thus we obtain, $p_i(\zeta_1,\zeta_2)=\zeta_1$ for $1\leq i\leq 8.$ Thus the $8^{th}$ approximate solution of Eq. (\ref{53}) is, \begin{equation} \label{56}
	\mathcal{V}_8(\zeta_1,\zeta_2,s)=\sum_{k=0}^8 \frac{\zeta_1}{s^{k\gamma+1}}.
\end{equation} Operating the inverse LT on the Eq. (\ref{56}). Then we get,  \begin{equation} \label{57}
\nu_8(\zeta_1,\zeta_2,\tau)= 
\sum_{k=0}^8 \frac{\zeta_1 \tau ^{k\gamma}}{\Gamma(k\gamma +1)}.
\end{equation} 
The infinite series solution of concerned problem is, \begin{equation} \label{59}
\nu(\zeta_1,\zeta_2,\tau)=\sum_{k=0}^\infty \frac{\zeta_1 \tau ^{k\gamma}}{\Gamma(k\gamma+1)} .
\end{equation} The closed-form solution is in the form $\nu (\zeta_1,\zeta_2, \tau)=\zeta_1 E_\gamma(\tau^\gamma).$ For $\gamma=1,$ the series reduces to $\nu(\zeta_1,\zeta_2,\tau)=\zeta_1\exp(\tau)$.
  
\end{solution}
\end{example}

\begin{example}   \label{ex6}
	Consider the 2-dimensional non-linear TFFPE, \begin{equation} \label{61}
		\begin{split}
		D_\tau^{\gamma}\nu(\zeta_1,\zeta_2,\tau)=&-D_{\zeta_1}\bigg(\frac{4}{\zeta_1}\nu^2(\zeta_1,\zeta_2,\tau)\bigg)-D_{\zeta_2}(\zeta_2\nu(\zeta_1,\zeta_2,\tau))+D_{\zeta_1\zeta_1}(\nu^2(\zeta_1,\zeta_2,\tau))+D_{\zeta_1\zeta_2}\nu(\zeta_1,\zeta_2,\tau)\\
		&+D_{\zeta_2\zeta_1}\nu(\zeta_1,\zeta_2,\tau)+D_{\zeta_2\zeta_2}\nu(\zeta_1,\zeta_2,\tau)	
	\end{split}
	\end{equation}
\begin{equation} \label{61a}
\nu(\zeta_1,\zeta_2,0)=\zeta_1^2, \hspace{0.4cm} (\zeta_1, \zeta_2) \in\mathbb{R}^2,\hspace{0.2cm} \tau \geq 0.
\end{equation}
The exact solution is $\nu(\zeta_1,\zeta_2,\tau)=\zeta_1^2E_\gamma(-\tau^{\gamma}).$\begin{solution}
    By using LT in Eq. (\ref{61}) and substituting (\ref{61a}), we obtain \begin{equation} \label{62}
\begin{split}
\mathcal{V}(\zeta_1,\zeta_2,s)=\ &\frac{\zeta_1^2}{s}-\frac{1}{s^\gamma}\bigg[ D_{\zeta_1}\biggr(\frac{4}{\zeta_1}\mathcal{L}\{({\mathcal{L}^{-1}\{\mathcal{V}(\zeta_1,\zeta_2,s)\})^2}\}\biggl)+D_{\zeta_2}(\zeta_2\mathcal{V}(\zeta_1,\zeta_2,s))-D_{\zeta_1\zeta_1}(\mathcal{L}\{(\mathcal{L}^{-1}\{\mathcal{V}(\zeta_1,\zeta_2,s)\})^2\})
\\&-D_{\zeta_1\zeta_2}\mathcal{V}(\zeta_1,\zeta_2,s)-D_{\zeta_2\zeta_1}\mathcal{V}(\zeta_1,\zeta_2,s)-D_{\zeta_2\zeta_2}\mathcal{V}(\zeta_1,\zeta_2,s)\bigg].\\
\end{split}
\end{equation} Repeating the procedure as in preceding example and 
 utilizing property \textbf{P4} for each $k=1,2,3\hdots$ we obtain $p_i(\zeta_1, \zeta_2)=(-1)^i\zeta_1^2$ for $1\leq i\leq 8.$ Thus, the $8^{th}$ approximate solution $\mathcal{V}_8$ of Eq. (\ref{62}) is presented as, \begin{equation}
\mathcal{V}_8(\zeta_1,\zeta_2,s)=
\sum_{k=0}^8 (-1)^k \frac{\zeta_1^2}{s^{k\gamma+1}} .
\end{equation} Applying inverse LT, we obtain \begin{equation}
\nu_8(\zeta_1,\zeta_2,\tau) 
= \sum_{k=0}^8 (-1)^k \frac{\zeta_1^2 \tau^{k\tau}}{\Gamma(k\gamma+1)}.
\end{equation} In general, the infinite series solution of Eq. (\ref{62}) have the form \begin{equation}
\mathcal{V}(\zeta_1,\zeta_2,s)=\sum_{k=0}^\infty (-1)^k \frac{\zeta_1^2 \tau^{k\tau}}{\Gamma(k\gamma+1)}.
\end{equation} 

Thus, after applying inverse LT, we have
  $\nu(\zeta_1,\zeta_2,\tau)$ $=\zeta_1^2E_\gamma(-\tau^\gamma).$ If $\gamma=1,$ this series reduces to $\nu(\zeta_1,\zeta_2,\tau)=\zeta_1^2\exp({-\tau}).$

\end{solution} 
\end{example}

\begin{example}   \label{ex7}
	Consider the 3-dimensional linear TFFPE, \begin{equation} \label{65}
		\begin{split}
		D_\tau^\gamma \nu(\zeta_1,\zeta_2,\zeta_3,\tau)=&-D_{\zeta_1}(2\zeta_1\nu(\zeta_1,\zeta_2,\zeta_3,\tau))-D_{\zeta_2}(2\zeta_2\nu(\zeta_1,\zeta_2,\zeta_3,\tau))-D_{\zeta_3}(2\zeta_3\nu(\zeta_1,\zeta_2,\zeta_3,\tau))\\
		&
		+D_{\zeta_1\zeta_1}(\zeta_1\nu(\zeta_1,\zeta_2,\zeta_3,\tau))+D_{\zeta_1\zeta_2}\nu(\zeta_1,\zeta_2,\zeta_3,\tau)+D_{\zeta_1\zeta_3}\nu(\zeta_1,\zeta_2,\zeta_3,\tau)\\
		&+D_{\zeta_2\zeta_1}\nu(\zeta_1,\zeta_2, \zeta_3,\tau)+D_{\zeta_2\zeta_2}\nu(\zeta_1,\zeta_2,\zeta_3,\tau)+D_{\zeta_2\zeta_3}\nu(\zeta_1,\zeta_2,\zeta_3,\tau)+D_{\zeta_3\zeta_1}\nu(\zeta_1,\zeta_2,\zeta_3,\tau)\\
	&+D_{\zeta_3\zeta_2}\nu(\zeta_1,\zeta_2,\zeta_3,\tau)+D_{\zeta_3\zeta_3}\bigg(\frac{3}{2}\zeta_3^3\nu(\zeta_1,\zeta_2,\zeta_3,\tau)\bigg),  \hspace{0.2cm} (\zeta_1, \zeta_2, \zeta_3) \in\mathbb{R}^3,\hspace{0.2cm} \tau \geq 0,
		\end{split}
	\end{equation} subject to the initial condition \begin{equation}
	\nu(\zeta_1,\zeta_2,\zeta_3,0)=\zeta_3.
\end{equation} Here $\gamma \in (0,1]$ and the exact solution is $\nu(\zeta_1,\zeta_2,\zeta_3,\tau)=\zeta_3 E_\gamma(\tau^\gamma).$\begin{solution}
     Operating the LT on Eq. (\ref{65}) 
to reduce the FDE into the following algebraic equation \begin{equation}
	\begin{split}
	\mathcal{V}(\zeta_1,\zeta_2,\zeta_3,\tau)=\ &\frac{\zeta_3}{s}-\frac{1}{s^\gamma}\Bigg[D_{\zeta_1}(2\zeta_1\mathcal{V}(\zeta_1,\zeta_2,\zeta_3,s))+D_{\zeta_2}(2\zeta_2\mathcal{V}(\zeta_1,\zeta_2,\zeta_3,s))+D_{\zeta_3}(2\zeta_3\mathcal{V}(\zeta_1,\zeta_2,\zeta_3,s))\\
	&-D_{\zeta_1\zeta_1}(\zeta_1\mathcal{V}(\zeta_1,\zeta_2,\zeta_3,s))-D_{\zeta_1\zeta_2}\mathcal{V}(\zeta_1,\zeta_2,\zeta_3,s)-D_{\zeta_1\zeta_3}\mathcal{V}(\zeta_1,\zeta_2,\zeta_3,s)-D_{\zeta_2\zeta_1}\mathcal{V}(\zeta_1,\zeta_2,\zeta_3,s)\\
	&-D_{\zeta_2\zeta_2}\mathcal{V}(\zeta_1,\zeta_2,\zeta_3,s)-D_{\zeta_2\zeta_3}\mathcal{V}(\zeta_1,\zeta_2,\zeta_3,s)-D_{\zeta_3\zeta_1}\mathcal{V}(\zeta_1,\zeta_2,\zeta_3,s)-D_{\zeta_3\zeta_2}\mathcal{V}(\zeta_1,\zeta_2,\zeta_3,s)\\
	&-D_{\zeta_3\zeta_3}\bigg(\frac{3}{2}\zeta_3^3\mathcal{V}(\zeta_1,\zeta_2,\zeta_3,s)\bigg)\Bigg].
		\end{split}
\end{equation} 
Now utilizing property \textbf{P4}, 
we obtain the unknown coefficient in the series representation as, $p_i(\zeta_1,\zeta_2,\zeta_3)=\zeta_3$ for $1\leqslant i \leqslant 8.$ Hence, the $8^{th}$ approximate solution of transformed equation is
\begin{equation} \label{68}
\mathcal{V}_8(\zeta_1,\zeta_2,\zeta_3,s)=\sum_{k=0}^8 \frac{\zeta_3}{s^{k\gamma +1}}.
\end{equation} Thus, the infinite series expression of the solution is 
\begin{equation}
\mathcal{V}(\zeta_1,\zeta_2,\zeta_3,s)=\sum_{k=0}^{\infty}\frac{\zeta_3}{s^{k\gamma+1}}.
\end{equation} Applying inverse LT in Eq. (\ref{68}), we get 
\begin{equation}
\nu(\zeta_1,\zeta_2,\zeta_3,\tau)=\zeta_3E_\gamma(\tau^\gamma).
\end{equation} For $\gamma=1,$ the solution converges to $\nu(\zeta_1,\zeta_2,\zeta_3,\tau)=\zeta_3\exp(\tau).$

\end{solution}
\end{example}

\begin{example}   \label{ex8}
	Consider the non-linear 3-dimensional TFFPE,\begin{equation}  \label{71}
		\begin{split}
				D_\tau^\gamma \nu(\zeta_1,\zeta_2,\zeta_3,\tau)=&-D_{\zeta_1}(\zeta_1\nu(\zeta_1,\zeta_2,\zeta_3,\tau))-D_{\zeta_2}(\nu^2(\zeta_1,\zeta_2,\zeta_3,\tau))-D_{\zeta_3}\bigg(\frac{2}{\zeta_3-1}\nu(\zeta_1,\zeta_2,\zeta_3,\tau)\bigg)\\
				&+D_{\zeta_1\zeta_1}(\nu^2(\zeta_1,\zeta_2,\zeta_3,\tau))+D_{\zeta_1\zeta_2}\nu(\zeta_1,\zeta_2,\zeta_3,\tau)+D_{\zeta_1\zeta_3}\nu(\zeta_1,\zeta_2,\zeta_3,\tau)\\
				&+D_{\zeta_2\zeta_1}\nu(\zeta_1,\zeta_2,\zeta_3,\tau)+D_{\zeta_2\zeta_2}(\nu^2(\zeta_1,\zeta_2,\zeta_3,\tau))+D_{\zeta_2\zeta_3}\nu(\zeta_1,\zeta_2,\zeta_3,\tau)+D_{\zeta_3\zeta_1}\nu(\zeta_1,\zeta_2,\zeta_3,\tau)\\
				&+D_{\zeta_3\zeta_2}\nu(\zeta_1,\zeta_2,\zeta_3,\tau)+D_{\zeta_3\zeta_3}\nu(\zeta_1,\zeta_2,\zeta_3,\tau),  \hspace{0.2cm} (\zeta_1, \zeta_2) \in\mathbb{R}^2,\hspace{0.2cm} \tau \geq 0,
		\end{split}
	\end{equation} with the initial condition $\nu(\zeta_1,\zeta_2,\zeta_3,0)=(\zeta_3-1)^2.$ The solution of Eq. (\ref{71}) is $\nu(\zeta_1,\zeta_2,\zeta_3,\tau)=(\zeta_3-1)^2E_\gamma(\tau^\gamma).$
\begin{solution}
    Applying LT on the Eq. (\ref{71}), we have 

\begin{equation}
	\begin{split}
	\mathcal{V}(\zeta_1,\zeta_2,\zeta_3,s)=\ &\frac{(\zeta_3-1)^2}{s}-\frac{1}{s^\gamma}\Bigg[D_{\zeta_1}(\zeta_1\mathcal{V}(\zeta_1,\zeta_2,\zeta_3,s))+	D_{\zeta_2}(\mathcal{L}\{(\mathcal{L}^{-1}\{\mathcal{V}(\zeta_1,\zeta_2,\zeta_3,s)\})^2\})\\
	&+D_{\zeta_3}\bigg(\frac{2}{\zeta_3-1}\mathcal{V}(\zeta_1,\zeta_2,\zeta_3,s)\bigg)- D_{\zeta_1\zeta_1}(\mathcal{L}\{(\mathcal{L}^{-1}\{\mathcal{V}(\zeta_1,\zeta_2,\zeta_3,s)\})^2\})- D_{\zeta_1\zeta_2}\mathcal{V}(\zeta_1,\zeta_2,\zeta_3,s)\\
	&-D_{\zeta_1\zeta_3}\mathcal{V}(\zeta_1,\zeta_2,\zeta_3,s)- D_{\zeta_2\zeta_1}\mathcal{V}(\zeta_1,\zeta_2,\zeta_3,s)- D_{\zeta_2\zeta_2}(\mathcal{L}\{(\mathcal{L}^{-1}\{\mathcal{V}(\zeta_1,\zeta_2,\zeta_3,s)\})^2\})\\
	&-D_{\zeta_2\zeta_3}\mathcal{V}(\zeta_1,\zeta_2,\zeta_3,s)-D_{\zeta_3\zeta_1}\mathcal{V}(\zeta_1,\zeta_2,\zeta_3,s)-
	D_{\zeta_3\zeta_2}\mathcal{V}(\zeta_1,\zeta_2,\zeta_3,s)-D_{\zeta_3\zeta_3}\mathcal{V}(\zeta_1,\zeta_2,\zeta_3,s) \Bigg].	
\end{split}
\end{equation} 
 Applying property \textbf{P4}, 
  we get the unknown coefficients $p_i(\zeta_1,\zeta_2,\zeta_3)=(\zeta_3-1)^2$ for $1\leqslant i\leqslant 8.$ Hence the $8^{th}$ approximate solution of transformed equation is \begin{equation}
\mathcal{V}_8(\zeta_1,\zeta_2,\zeta_3,s)=\sum_{k=0}^8 \frac{(\zeta_3-1)^2}{s^{k\gamma+1}}.
\end{equation} Thus, the closed-form solution of the transformed equation is \begin{equation}
\mathcal{V}(\zeta_1,\zeta_2,\zeta_3,s)=\sum_{k=0}^{\infty}\frac{(\zeta_3-1)^2}{s^{k\gamma+1}}.
\end{equation} Finally, we apply inverse LT to obtain the desired solution 
\begin{equation} \label{67}
	\nu(\zeta_1,\zeta_2,\zeta_3,\tau)=(\zeta_3-1)^2E_\gamma(\tau^\gamma).
\end{equation} For $\gamma=1,$ Eq. (\ref{67}) reduces to the exact solution $\nu(\zeta_1,\zeta_2,\zeta_3,\tau)=(\zeta_3-1)^2\exp (\tau).$

\end{solution} 
\end{example}

\section{Discussions}

\begin{figure}
	\centering
	\subfloat[$\nu_8(\zeta,\tau),\ \gamma=0.5$.]{%
		\resizebox*{6cm}{!}{\includegraphics{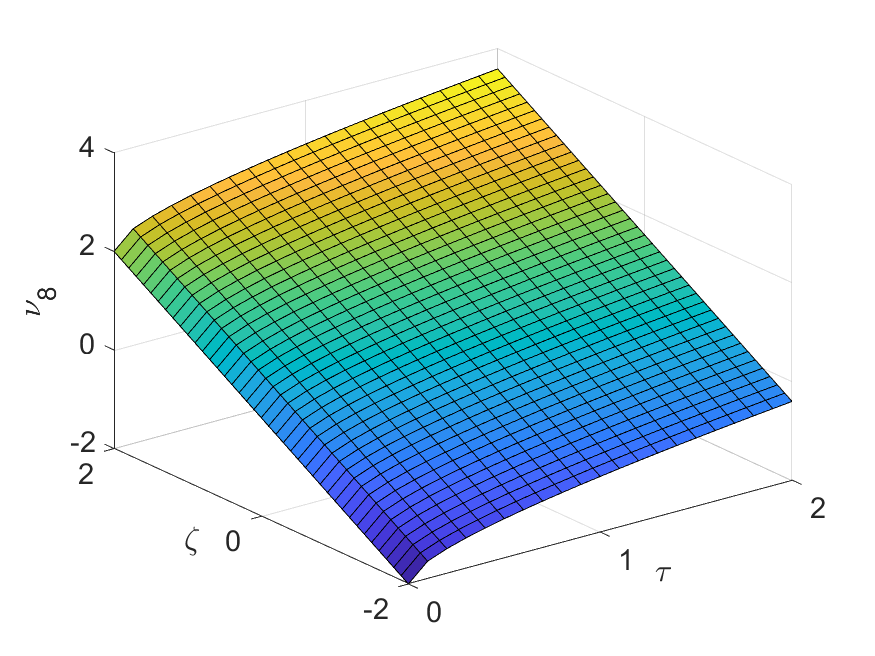}}}
	\subfloat[$\nu_8(\zeta,\tau),\ \gamma=0.75$.]{%
		\resizebox*{6cm}{!}{\includegraphics{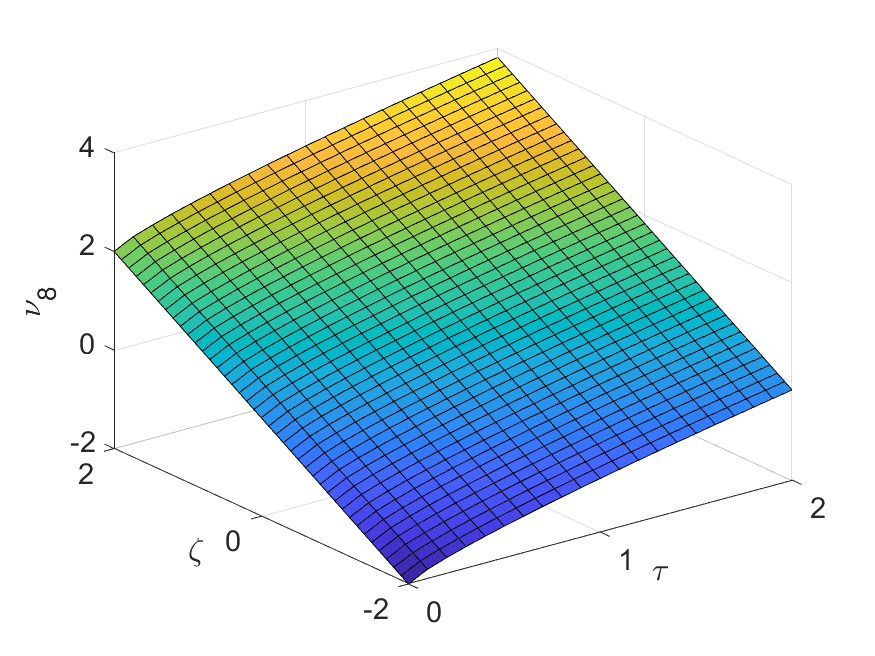}}}\hspace{5pt}
	\subfloat[$\nu_8(\zeta,\tau),\ \gamma=1$.]{%
		\resizebox*{6cm}{!}{\includegraphics{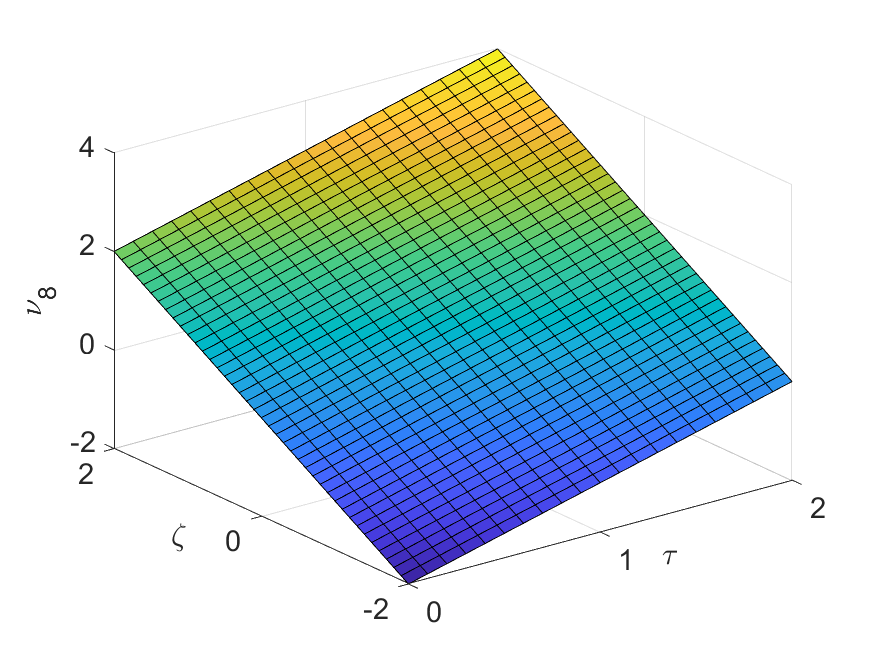}}}
	\subfloat[$\nu(\zeta,\tau),\  \gamma=1$.]{%
		\resizebox*{6cm}{!}{\includegraphics{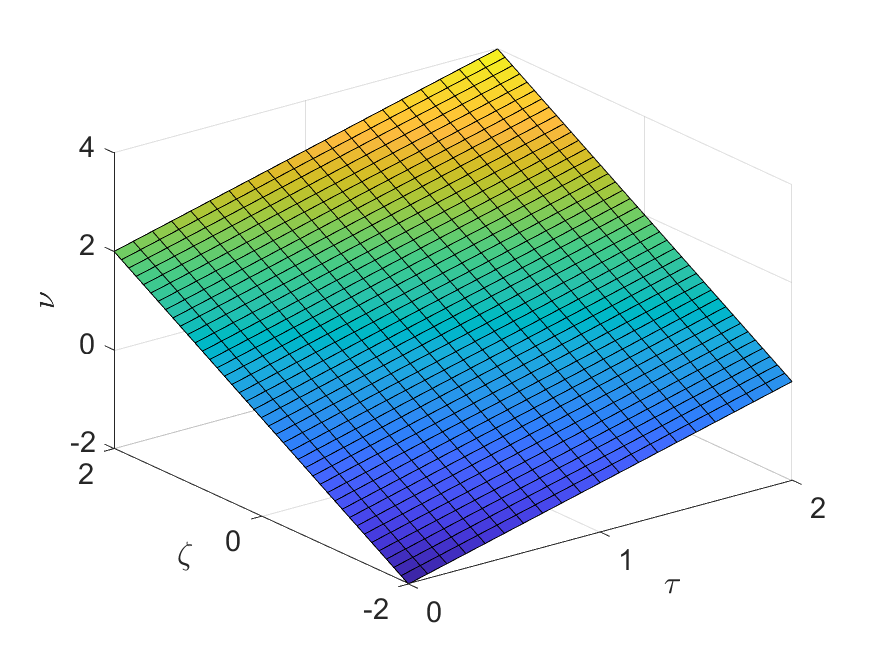}}}
	\caption{The variation of solution profile for various fractional orders for Example \ref{ex1}.} \label{fig_1}
\end{figure}
\begin{figure}
	\centering
	\includegraphics[width=9cm]{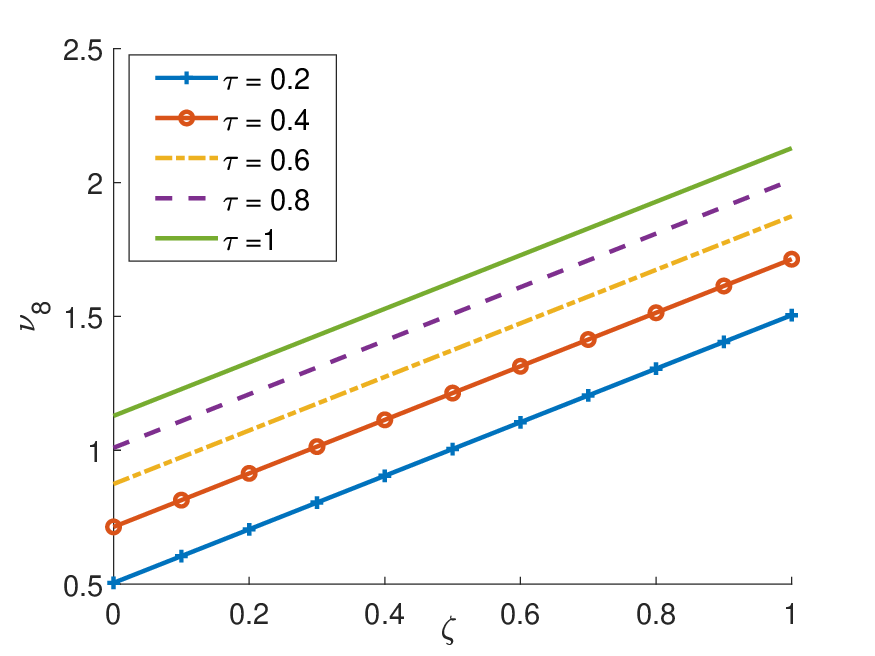}
	\caption{The variation of $\nu_8(\zeta,\tau)$ for different time $\tau$ for Example \ref{ex1}.} \label{fig_1a}
\end{figure}


\begin{table}

\centering
\resizebox{16cm}{!}{
	{\begin{tabular}{cccccc}
			\toprule
			$\tau$ &$\gamma$ = 0.2&  $\gamma=0.4$ &  $\gamma=0.6$ & $\gamma=0.8$& $\gamma=1$\\
			\midrule 
			0.15 & 1.600081043809555 & 0.952981800749650 & 0.737755085043236 & 0.636417247831176 & 0.580917121364088 \\
			0.30 & 2.135367406482500 & 1.247856716099542 & 0.925069682147884 & 0.766730117890084 & 0.674929403760045 \\
			0.45 & 2.631768336071708 & 1.556097169766161 & 1.126536612111816 & 0.911564290106330 & 0.784156091653569 \\
			0.60 & 3.116845586716595 & 1.893695854737707 & 1.352081318394315 & 1.076288836236974 & 0.911059385428571 \\
			0.75 & 3.599855147026588 & 2.268773048642961 & 1.608451623447892 & 1.265354801929119 & 1.058499896526337 \\
			0.90 & 4.084939830475697 & 2.687165687471149 & 1.901880218811209 & 1.483342806808060 & 1.229800969474331 \\
			\bottomrule			
			
	\end{tabular}}
        }
\caption{LRPS solution $\nu_8(0.5,\tau)$ for distinct $\gamma$ values for Example \ref{ex1}. }	
	\label{tab:1}
\end{table}

\begin{table}
\centering
	{\begin{tabular}{cccc}
			\toprule
			$\zeta$& LRPS & AHPM \cite{HAS24}& Exact \\
			\midrule 
			
			-1.0     & -0.99900 & -0.99900 & -0.99900 \\
			-0.8     & -0.79900 & -0.79900 & -0.79900  \\
			-0.6     & -0.59900 & -0.59900 & -0.59900  \\
			-0.4     & -0.39900 & -0.39900 & -0.39900  \\
			-0.2     & -0.19900 & -0.19900 & -0.19900  \\
			0.0      &  0.00100 &  0.00100 &  0.00100  \\
			0.2      &  0.20100 &  0.20100 &  0.20100  \\
			0.4      &  0.40100 &  0.40100 &  0.40100  \\
			0.6      &  0.60100 &  0.60100 &  0.60100  \\
			0.8      &  0.80100 &  0.80100 &  0.80100  \\
			1.0      &  1.00100 &  1.00100 &  1.00100  \\
			\bottomrule			
			
	\end{tabular}}
	\caption{Solution comparisons 
		at $\tau=0.001$ and $\gamma=1$ for Example \ref{ex1}.}
	\label{tab:2}
\end{table}


\begin{figure}
	\centering
	
	\subfloat[$\nu_8(\zeta,\tau),\ \gamma=0.5$.]{%
		\resizebox*{6cm}{!}{\includegraphics{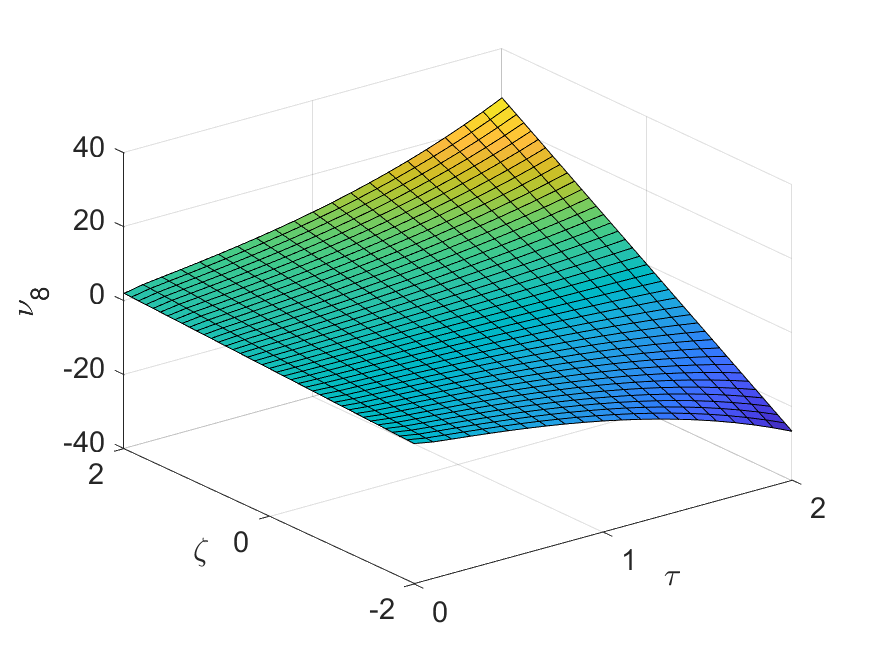}}}
	\subfloat[$\nu_8(\zeta,\tau),\ \gamma=0.75$.]{%
		\resizebox*{6cm}{!}{\includegraphics{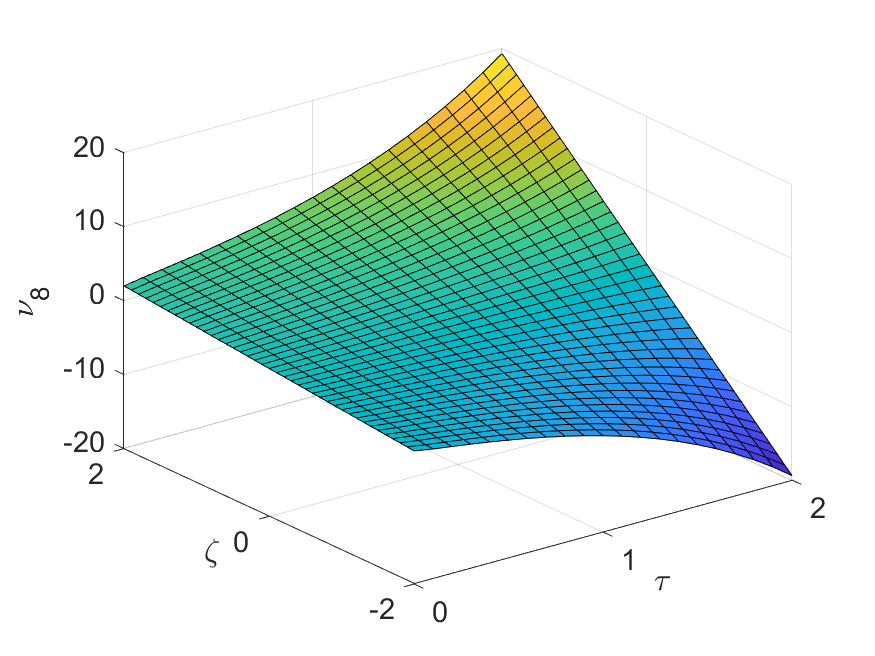}}}\hspace{5pt}
	\subfloat[$\nu_8(\zeta,\tau),\ \gamma=1$.]{%
		\resizebox*{6cm}{!}{\includegraphics{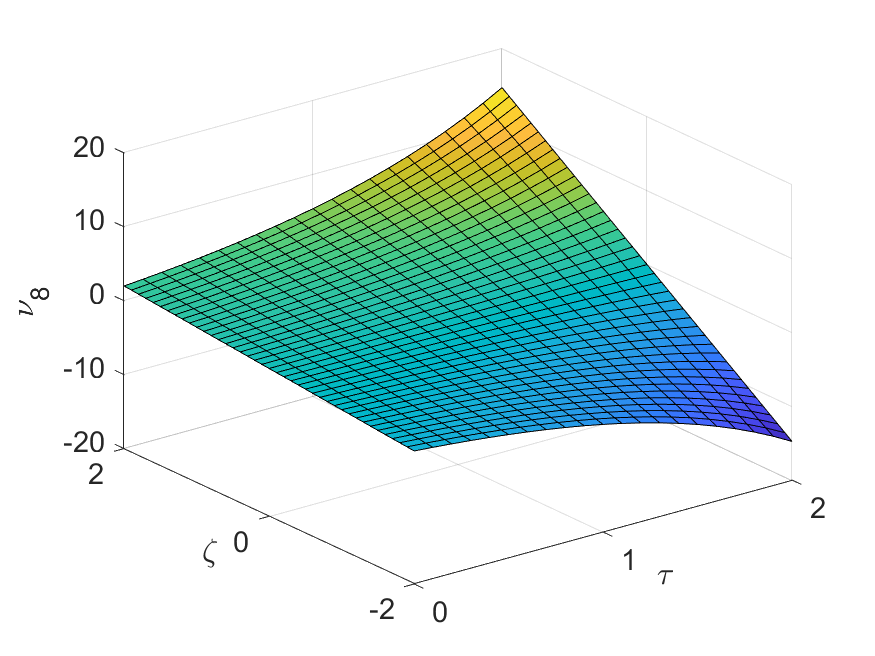}}}
	\subfloat[$\nu(\zeta,\tau),\  \gamma=1.$]{%
		\resizebox*{6cm}{!}{\includegraphics{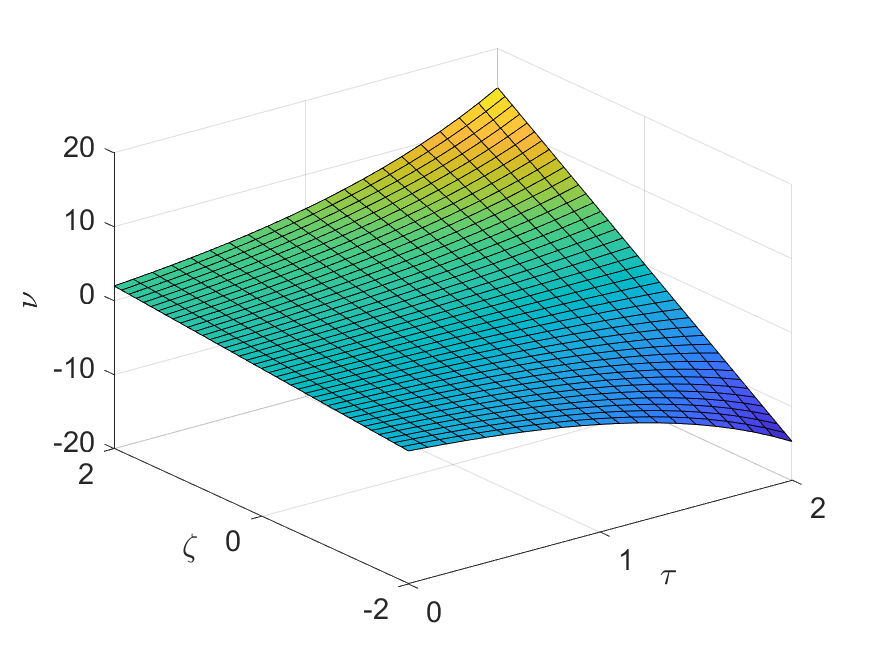}}}
	\caption{The variation of solution profile for various fractional orders for Example \ref{ex2}.} \label{fig_2}
\end{figure}
	\begin{figure}
	\centering
	\includegraphics[width=9cm]{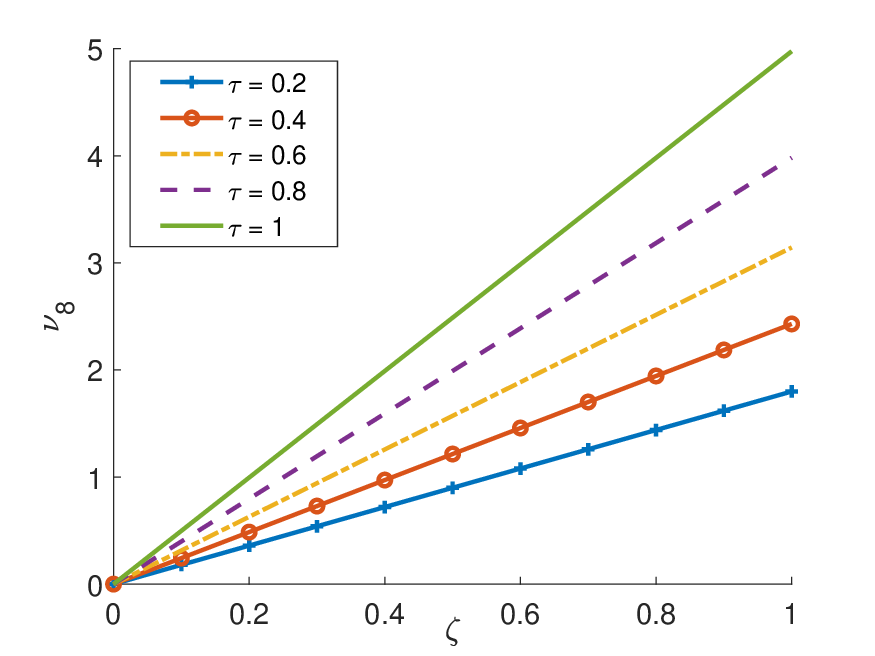}
	\caption{The variation of $\nu_8(\zeta,\tau)$ for different time $\tau$ for Example \ref{ex2}.} \label{fig_2a}
	\end{figure}

\begin{figure}
	\centering
	
	\subfloat[$|\nu(\zeta,\tau)-\nu_8(\zeta,\tau)|$ for $\gamma=0.5.$]{%
		\resizebox*{6cm}{!}{\includegraphics{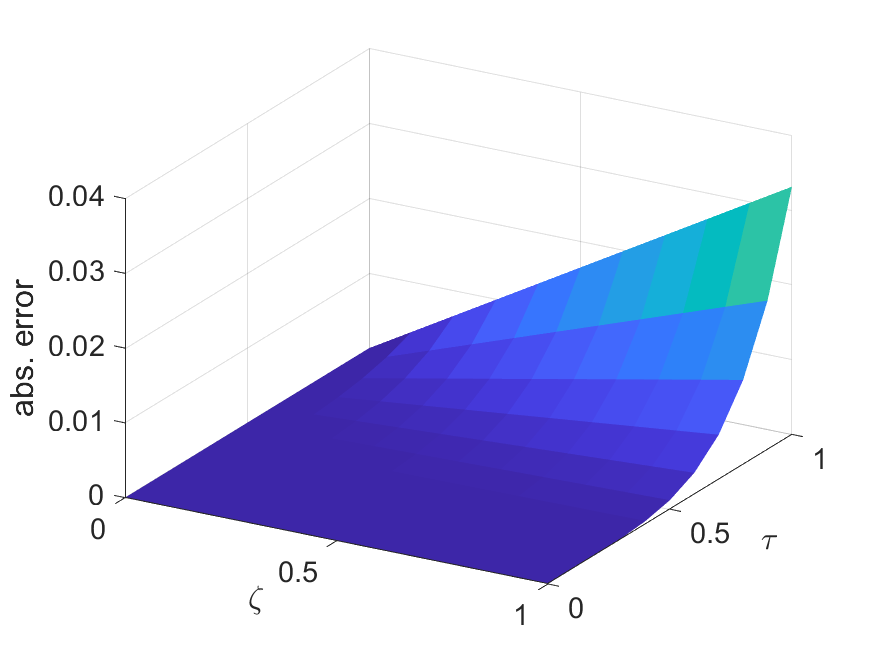}}}\hspace{1cm}
	\subfloat[$\zeta=0.5, \ \gamma=1$.]{%
		\resizebox*{6cm}{!}{\includegraphics{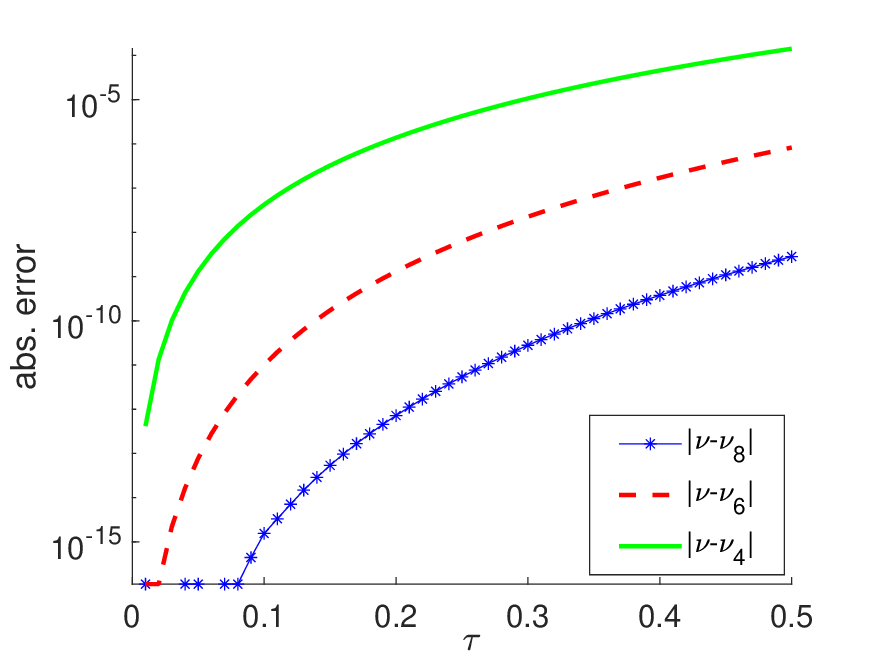}}}
	\caption{Error plots for Example \ref{ex2}.} 
    \label{fig_3}
\end{figure}


\begin{table}
\centering
    \resizebox{16cm}{!}{
	{\begin{tabular}{cccccc}
			\toprule
			$\tau$ & $\gamma$ = 0.2& $\gamma=0.4$ &  $\gamma=0.6$ & $\gamma=0.8$& $\gamma=1$\\
			\midrule
			0.15 & 1.600081043809555 & 0.952981800749650 & 0.737755085043236 & 0.636417247831176 & 0.580917121364088 \\
			0.30 & 2.135367406482500 & 1.247856716099542 & 0.925069682147884 & 0.766730117890084 & 0.674929403760045 \\
			0.45 & 2.631768336071708 & 1.556097169766161 & 1.126536612111816 & 0.911564290106330 & 0.784156091653569 \\
			0.60 & 3.116845586716595 & 1.893695854737707 & 1.352081318394315 & 1.076288836236974 & 0.911059385428571 \\
			0.75 & 3.599855147026588 & 2.268773048642961 & 1.608451623447892 & 1.265354801929119 & 1.058499896526337 \\
			0.90 & 4.084939830475697 & 2.687165687471149 & 1.901880218811209 & 1.483342806808060 & 1.229800969474331 \\
			\bottomrule			
			
	\end{tabular}}
    }
	\caption{LRPS solution $\nu_8(0.5,\tau)$ for various values of $\gamma$ for Example \ref{ex2}. }
	\label{tab:3}
\end{table}


\begin{table}
	\centering
    
    \resizebox{16cm}{!}{
		\begin{tabular}{ccccccc}
			\toprule
			$\zeta$ & $\tau$ & RPSM \cite{YKK15} & FRDTM \cite{DC15} & HPTM \cite{K13}&   LRPS & Exact \\
			\midrule
			0.25 & 0.01 & 0.2525125417 & 0.2525125417  & 0.2525125417 & 0.2525125417 & 0.2525125418 \\
			0.50 & & 0.5050250833 & 0.5050250833 & 0.5050250833 & 0.5050250833 & 0.5050250835 \\
			0.75 &  & 0.7575376250 & 0.7575376250 & 0.7575376250 & 0.7575376250 & 0.7575376252 \\
            
			1.00 & & 1.010050167 & 1.010050167 & 1.010050167 & 1.010050167 &  1.010050167 \\
            \midrule
			0.25 & 0.6 & 0.45400000000 & 0.45400000000 & 0.45400000000 & 0.45400000000 & 0.455297000 \\
			0.50 & & 0.90800000000 & 0.90800000000 & 0.90800000000 & 0.90800000000 & 0.911059400 \\
			0.75 & & 1.36200000000 & 1.36200000000 &  1.36200000000 & 1.36200000000 & 1.366589100 \\
			1.00 & & 1.81600000000 & 1.81600000000 & 1.81600000000 &  1.81600000000 & 1.822118800 \\
			\bottomrule
	\end{tabular}}
    \caption{Comparison of LRPS method with other existing methods for $k=3$ and $\gamma=1$ for Example \ref{ex2}.}
	\label{tab:4}
\end{table}


\begin{table}
\centering
	{\begin{tabular}{ccccc}
			\toprule
			$\tau$ &exact & LRPS & abs. error & relative error \\
			\midrule
			0.15 & 0.580917121364088 & 0.580917121364142 & 5.373e-14 & 1.000e-14 \\
			0.30 & 0.674929403760045 & 0.674929403788002 & 2.795e-11 & 4.224e-11 \\
			0.45 & 0.784156091653569 & 0.784156092745084 & 1.091e-09 & 1.127e-09 \\
			0.60 & 0.911059385428571 & 0.911059400195254 & 1.476e-08 & 1.634e-08 \\
			0.75 & 1.058499896526337 & 1.058500008306337 & 1.117e-07 & 1.056e-07 \\
			0.90 & 1.229800969474331 & 1.229801555578475 & 5.861e-07 & 5.747e-07 \\
			\bottomrule

	\end{tabular}}
	\caption{
    Error estimates between $\nu(0.5, \ \tau)$ and $\nu_8(0.5,\ \tau)$ 
for $\gamma=1$ for Example \ref{ex2}.}
	\label{tab:5}
\end{table}

\begin{table}
\centering
	{\begin{tabular}{cccccc}
			\toprule
			$\zeta$&	$\tau$ &$|\nu(\zeta, \ \tau)-\nu_4(\zeta,\ \tau)|$ & $|\nu(\zeta,\tau)-\nu_6(\zeta,\ \tau)|$ & $|\nu(\zeta,\ \tau)-\nu_8(\zeta,\tau)|$ \\
			\midrule
			
			0.5&	 0.1 & 4.237e-08 & 1.005e-11 & 1.554e-15 \\
			&	0.2 & 1.379e-06 & 1.302e-09& 7.199e-13 \\
			&	0.3 & 1.065e-05 & 2.254e-08 & 2.796e-11 \\
			&	0.4 & 4.568e-05 & 1.710e-07 & 3.762e-10 \\
			&	0.5 & 1.419e-04 & 8.263e-07 & 2.832e-09 \\
			\midrule
			1&	0.1 & 8.474e-08 & 2.009e-11 & 3.109e-15 \\
			&	0.2 & 2.758e-06 & 2.605e-09 & 1.440e-12 \\
			&	0.3 & 2.131e-05 & 4.508e-08 & 5.591e-11 \\
			&	0.4 & 9.136e-05 & 3.421e-07 & 7.524e-10 \\
			&	0.5 & 2.838e-04 & 1.653e-06 & 5.664e-09 \\
			
			\bottomrule
	\end{tabular}}
	\caption{Comparison of absolute errors for Example \ref{ex2}.}
	\label{tab:6}
\end{table}

\begin{figure}
	\centering
	
	\subfloat[$\nu_8(\zeta, \tau),\ \gamma=0.5$.]{%
		\resizebox*{6cm}{!}{\includegraphics{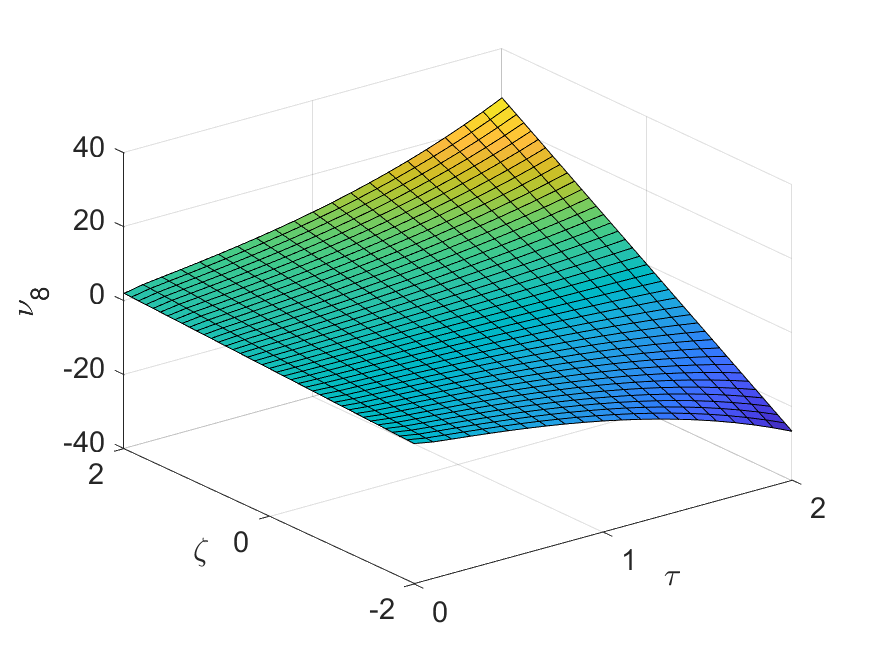}}}
	\subfloat[$\nu_8(\zeta,\tau),\ \gamma=0.75$.]{%
		\resizebox*{6cm}{!}{\includegraphics{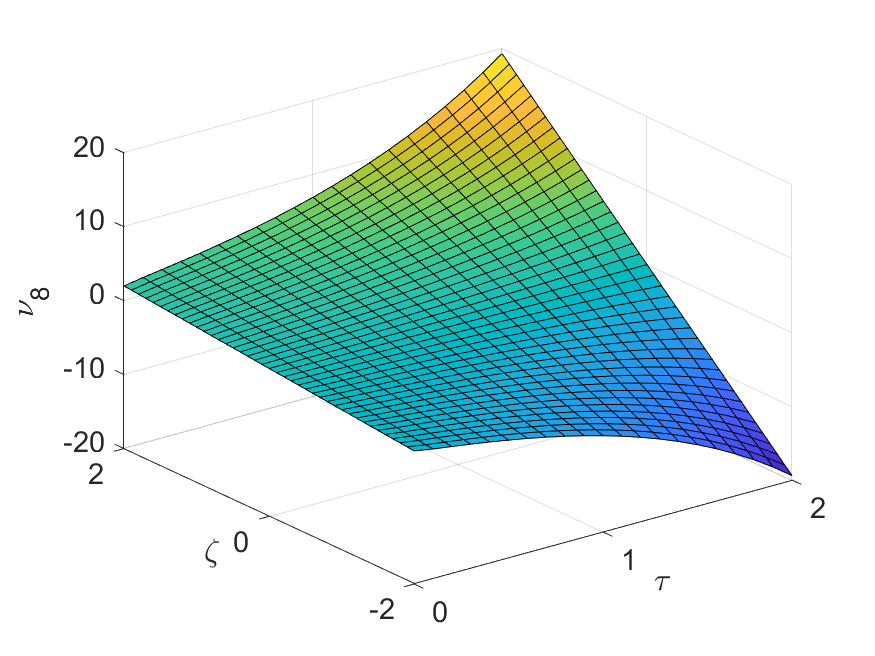}}}\hspace{5pt}
	\subfloat[$\nu_8(\zeta,\tau),\ \gamma=1$.]{%
		\resizebox*{6cm}{!}{\includegraphics{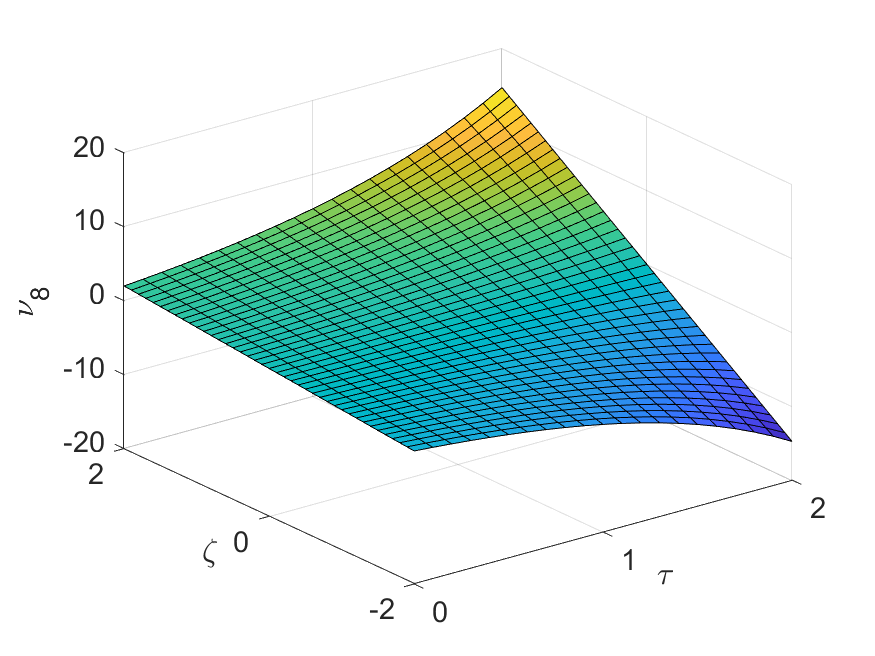}}}
	\subfloat[$\nu(\zeta,\tau),\ \gamma=1$.]{%
		\resizebox*{6cm}{!}{\includegraphics{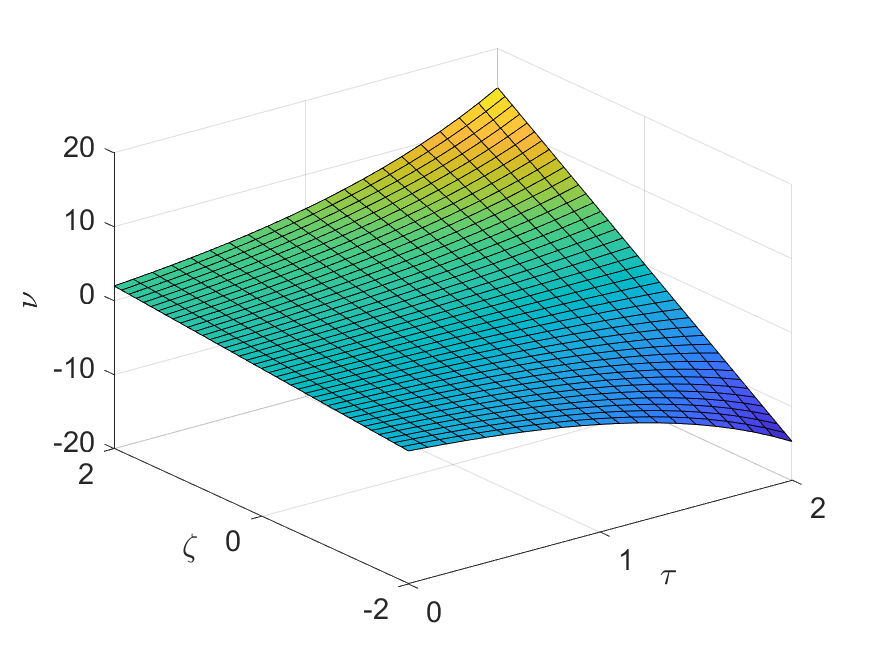}}}
	\caption{The variation of solution profile for various fractional orders for Example \ref{ex4}.} \label{fig_6}
\end{figure}

\begin{figure}
	\centering
    \includegraphics[width=9cm]{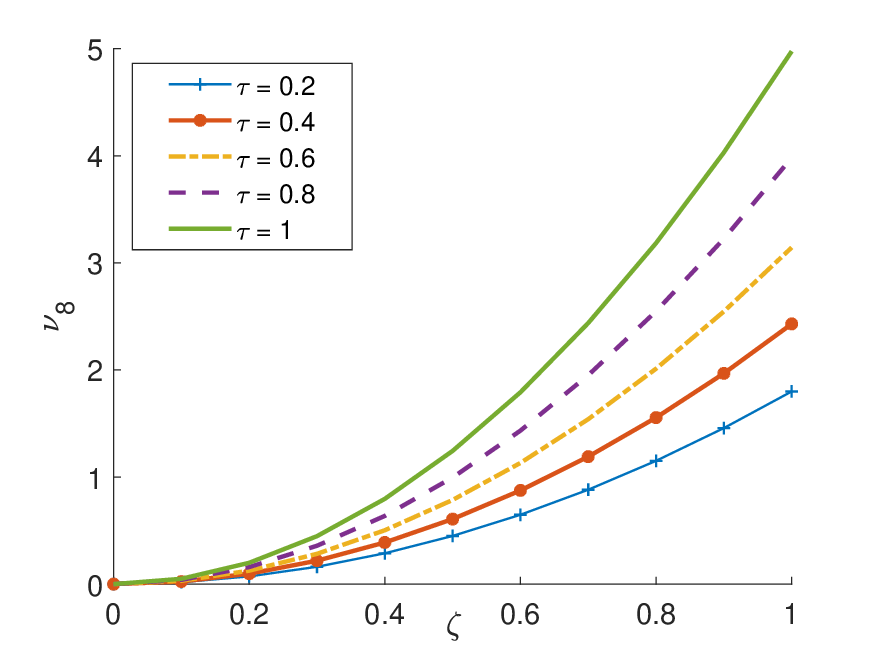}	
	\caption{The variation of $\nu_8(\zeta, \tau)$ for different time $\tau$ for Example \ref{ex4}.} \label{fig_4a}
\end{figure}
\begin{figure}
	\centering
	\subfloat[$|\nu(\zeta,\tau)-\nu_8(\zeta,\tau)|$ for $\gamma=0.5.$]{%
		\resizebox*{6cm}{!}{\includegraphics{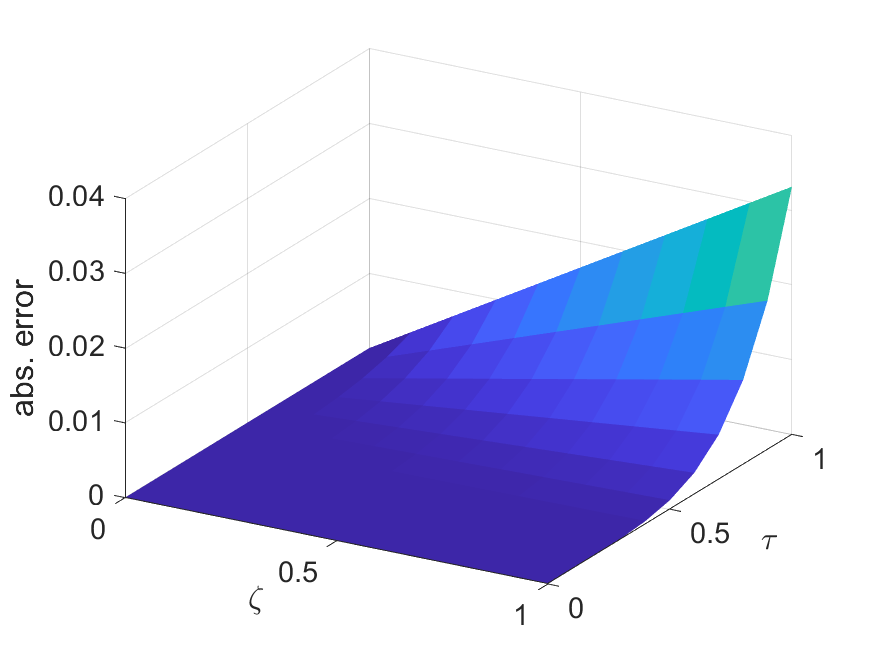}}} \hspace{1cm}
	\subfloat[$\zeta=0.5, \ \gamma=1$.]{%
		\resizebox*{6cm}{!}{\includegraphics{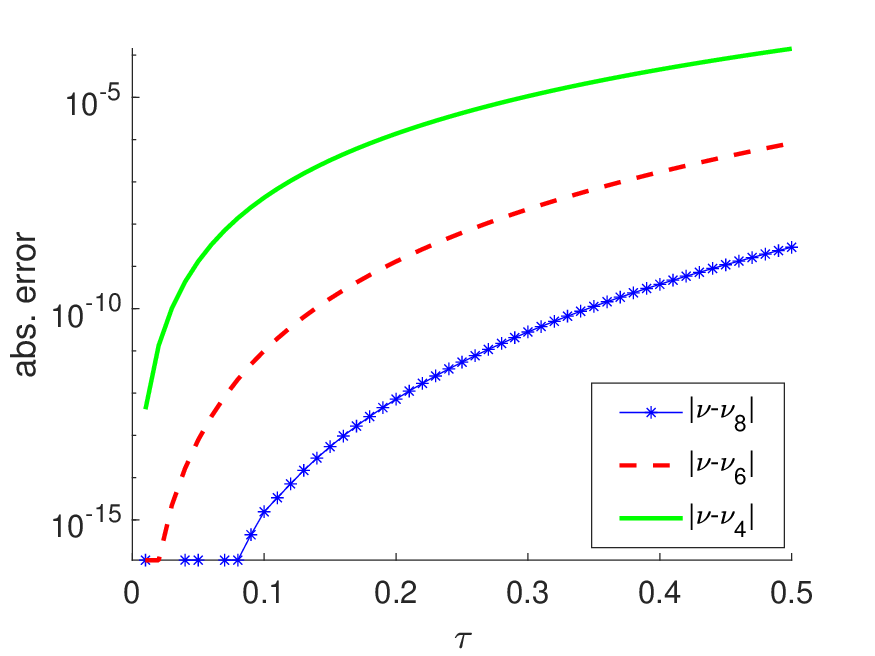}}}
	\caption{Error plots for Example \ref{ex4}.} \label{fig_7}
\end{figure}



\begin{table}
\centering
 \resizebox{16cm}{!}{
	
	\begin{tabular}{cccccc}
			\toprule
			$\tau$ & $\gamma=0.2$ & $\gamma=0.4$ & $\gamma=0.6$ & $\gamma=0.8$ & $\gamma=1$ \\
			\midrule
			0.15 & 1.600081043809555 & 0.952981800749650 & 0.737755085043236 & 0.636417247831176 & 0.580917121364088 \\
			0.30 & 2.135367406482500 & 1.247856716099542 & 0.925069682147884 & 0.766730117890084 & 0.674929403760045 \\
			0.45 & 2.631768336071708 & 1.556097169766161 & 1.126536612111816 & 0.911564290106330 & 0.784156091653569 \\
			0.60 & 3.116845586716595 & 1.893695854737707 & 1.352081318394315 & 1.076288836236974 & 0.911059385428571 \\
			0.75 & 3.599855147026588 & 2.268773048642961 & 1.608451623447892 & 1.265354801929119 & 1.058499896526337 \\
			0.90 & 4.084939830475697 & 2.687165687471149 & 1.901880218811209 & 1.483342806808060 & 1.229800969474331 \\
			\bottomrule
			
	\end{tabular}}
	\caption{LRPS solution $\nu_8(0.5,\ \tau)$ for distinct $\gamma$ values for Example \ref{ex4}.}
	\label{tab:11}
\end{table}



\begin{table}
\centering
	{\begin{tabular}{ccccc}
			\toprule
			$\tau$ &exact & LRPS & abs. error & relative error \\
			\midrule
			0.15 & 0.580917121364088 & 0.580917121364142 & 5.373e-14 & 1.000e-14 \\
			0.30 & 0.674929403760045 & 0.674929403788002 & 2.795e-11 & 4.224e-11 \\
			0.45 & 0.784156091653569 & 0.784156092745084 & 1.091e-09 & 1.127e-09 \\
			0.60 & 0.911059385428571 & 0.911059400195254 & 1.476e-08 & 1.634e-08 \\
			0.75 & 1.058499896526337 & 1.058500008306337 & 1.117e-07 & 1.056e-07 \\
			0.90 & 1.229800969474331 & 1.229801555578475 & 5.861e-07 & 5.747e-07 \\
			\bottomrule

	\end{tabular}}
         \caption{
          Error estimates between $\nu(0.5, \ \tau)$ and $\nu_8(0.5,\ \tau)$ for $\gamma=1$ for Example \ref{ex4}.}
	\label{tab:12}
\end{table}

\begin{table}
\centering
	{\begin{tabular}{cccccc}
			\toprule
			$\zeta$&	$\tau$ &$|\nu(\zeta,\ \tau)-\nu_4(\zeta,\ \tau)|$ & $|\nu(\zeta,\ \tau)-\nu_6(\zeta,\ \tau)|$ & $|\nu(\zeta,\ \tau)-\nu_8(\zeta,\ \tau)|$ \\
			\midrule
			
			0.5&	 0.1 & 4.237e-08 & 1.005e-11 & 1.554e-15 \\
			&	0.2 & 1.379e-06 & 1.302e-09 & 7.199e-13 \\
			&	0.3 & 1.065e-05 &2.254e-08 & 2.796e-11 \\
			&	0.4 & 4.568e-05 & 1.710e-07 & 3.762e-10 \\
			&	0.5 & 1.419e-04 & 8.263e-07 & 2.832e-09 \\
			\midrule
			1&	0.1 & 8.474e-08 & 2.009e-11 & 3.109e-15 \\
			&	0.2 & 2.758e-06 & 2.605e-09 & 1.440e-12 \\
			&	0.3 & 2.131e-05 & 4.508e-08 & 5.591e-11 \\
			&	0.4 & 9.136e-05 &3.421e-07 & 7.524e-10 \\
			&	0.5 & 2.838e-04 & 1.653e-06 & 5.664e-09 \\
			
			\bottomrule
	\end{tabular}}
	\caption{ Comparison of absolute errors for Example \ref{ex4}.}
	\label{tab:13}
\end{table}


\begin{figure}
	\centering
	\subfloat[$\nu_8(\zeta_1,\ \zeta_2,\ 0.5),\ \gamma=0.5$.]{%
		\resizebox*{6cm}{!}{\includegraphics{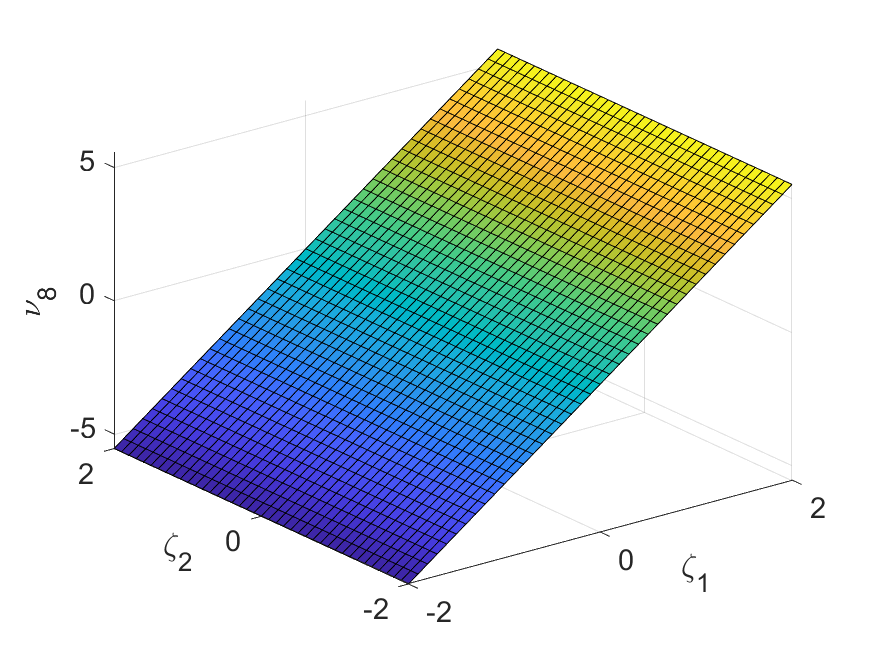}}}
	\subfloat[$\nu_8(\zeta_1,\ \zeta_2,\ 0.5),\ \gamma=0.75$.]{%
		\resizebox*{6cm}{!}{\includegraphics{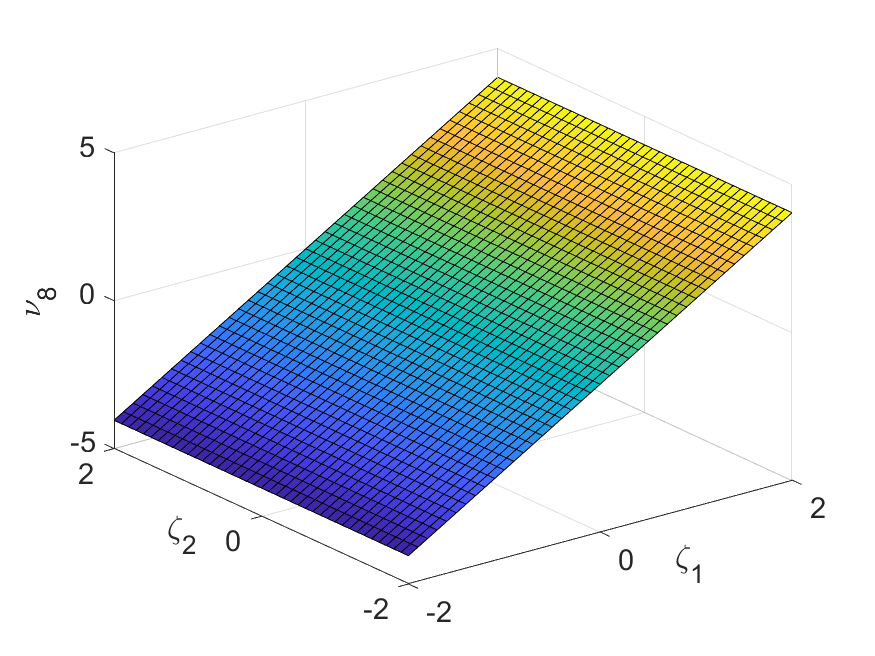}}}\hspace{5pt}
	\subfloat[$\nu_8(\zeta_1,\ \zeta_2,\ 0.5),\ \gamma=1$.]{%
		\resizebox*{6cm}{!}{\includegraphics{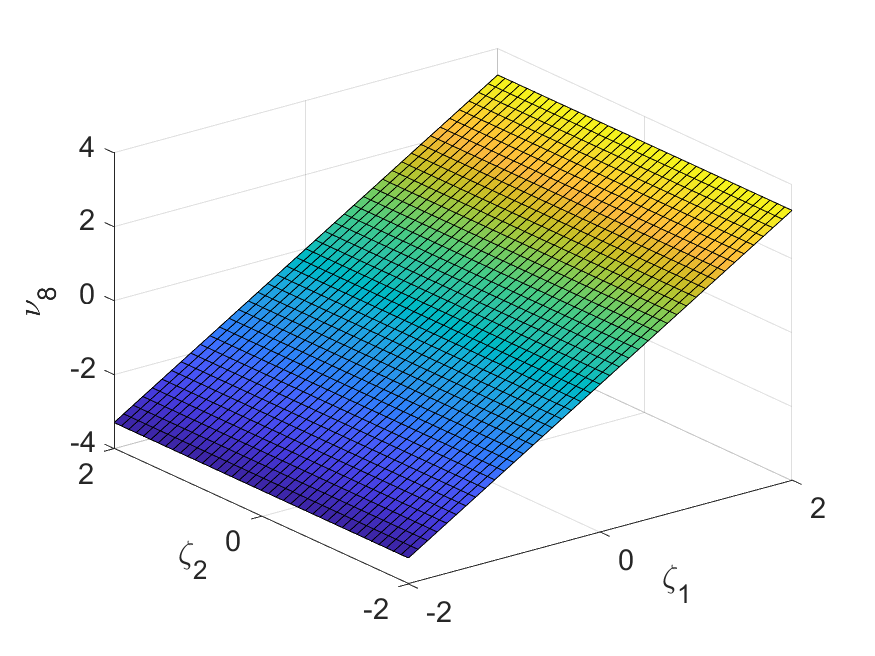}}}
	\subfloat[$\nu(\zeta_1,\ \zeta_2,\ 0.5),\ \gamma=1$.]{%
		\resizebox*{6cm}{!}{\includegraphics{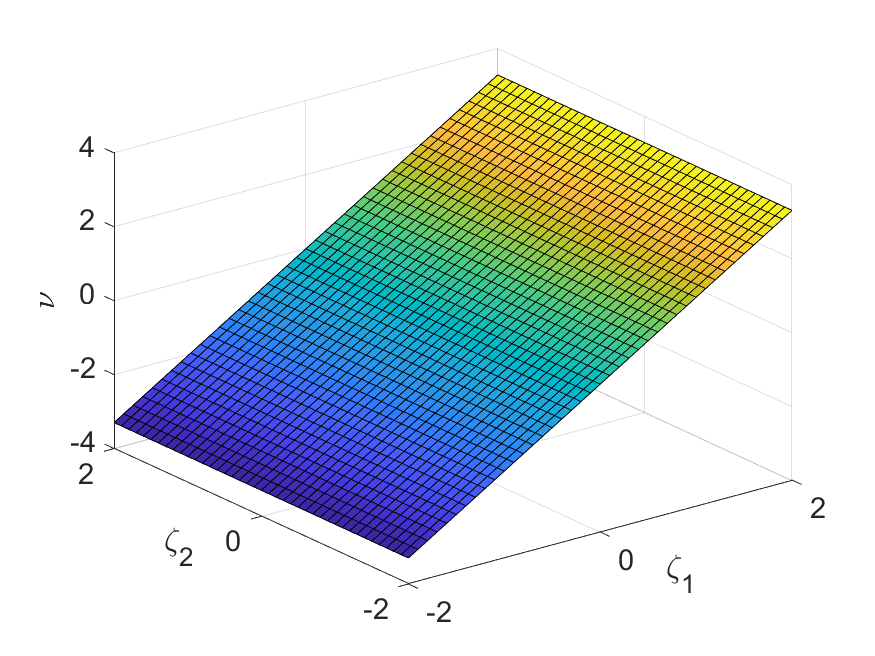}}}
	\caption{The variation of solution profile for various fractional orders for Example \ref{ex5}.} \label{fig_8}
\end{figure}

\begin{figure}
	\centering
	
	\subfloat[$|\nu(\zeta_1, \ \zeta_2, \ 0.5)-\nu_8(\zeta_1, \ \zeta_2, \ 0.5)|$ for $\gamma=0.5.$ ]{%
		\resizebox*{6cm}{!}{\includegraphics{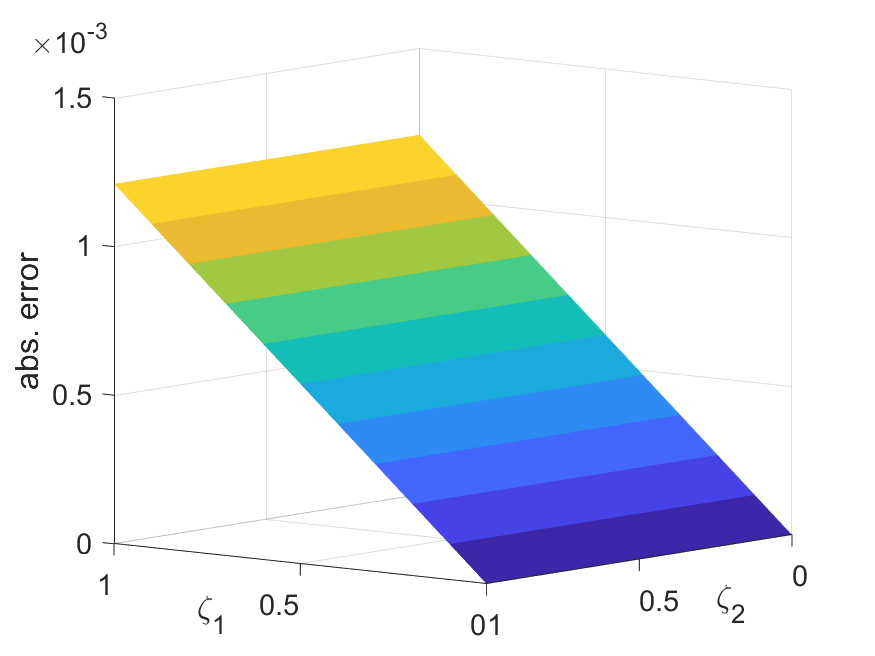}}} \hspace{1cm}
	\subfloat[$\zeta_1=\zeta_2=0.5,\ \gamma=1.$]{%
		\resizebox*{6cm}{!}{\includegraphics{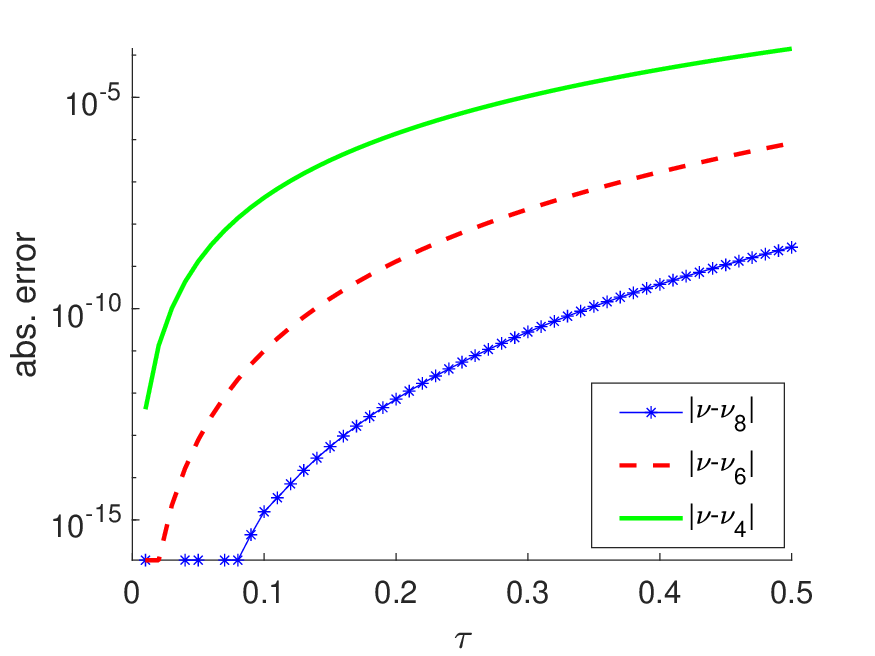}}}
	\caption{Error plots for Example \ref{ex5}.} \label{fig_9}
\end{figure}


\begin{table}
\centering
 \resizebox{16cm}{!}{
	\begin{tabular}{cccccc}
			\toprule
			$\tau$ & $\gamma=0.2$ & $\gamma=0.4$ & $\gamma=0.6$ & $\gamma=0.8$ & $\gamma=1$ \\
			\midrule
			0.15 & 1.600081043809555 & 0.952981800749650 & 0.737755085043236 & 0.636417247831176 & 0.580917121364088 \\
			0.30 & 2.135367406482500 & 1.247856716099542 & 0.925069682147884 & 0.766730117890084 & 0.674929403760045 \\
			0.45 & 2.631768336071708 & 1.556097169766161 & 1.126536612111816 & 0.911564290106330 & 0.784156091653569 \\
			0.60 & 3.116845586716595 & 1.893695854737707 & 1.352081318394315 & 1.076288836236974 & 0.911059385428571 \\
			0.75 & 3.599855147026588 & 2.268773048642961 & 1.608451623447892 & 1.265354801929119 & 1.058499896526337 \\
			0.90 & 4.084939830475697 & 2.687165687471149 & 1.901880218811209 & 1.483342806808060 & 1.229800969474331 \\
			\bottomrule
			
	\end{tabular}}
	\caption{LRPS solution $\nu_8(0.5,\ 0.5,\ \tau)$ for distinct $\gamma$ values for Example \ref{ex5}.}
	\label{tab:14}
\end{table}

\begin{table}
\centering
	{\begin{tabular}{ccccc}
			\toprule
			$\tau$ &exact & LRPS & abs. error & relative error \\
			\midrule
			0.15 & 0.580917121364088 & 0.580917121364142 & 5.373e-14 & 1.000e-14 \\
			0.30 & 0.674929403760045 & 0.674929403788002 & 2.795e-11 & 4.224e-11 \\
			0.45 & 0.784156091653569 & 0.784156092745084 & 1.091e-09 & 1.127e-09 \\
			0.60 & 0.911059385428571 & 0.911059400195254 & 1.476e-08 & 1.634e-08 \\
			0.75 & 1.058499896526337 & 1.058500008306337 & 1.117e-07 & 1.056e-07 \\
			0.90 & 1.229800969474331 & 1.229801555578475 & 5.861e-07 & 5.747e-07 \\
			\bottomrule

	\end{tabular}}
	\caption{
    Error estimates between $\nu(0.5, \ 0.5,\ \tau)$ and $\nu_8(0.5,\ 0.5,\ \tau)$ for $\gamma=1$ for Example \ref{ex5}.}
	\label{tab:15}
\end{table}

\begin{table}
\centering
	{\begin{tabular}{cccccc}
			\toprule
			$\zeta_1=\zeta_2$&	$\tau$ &$|\nu(\zeta_1,\ \zeta_2,\ \tau)-\nu_4(\zeta_1,\ \zeta_2,\ \tau)|$ & $|\nu(\zeta_1,\ \zeta_2,\ \tau)-\nu_6(\zeta_1,\ \zeta_2,\ \tau)|$ & $|\nu(\zeta_1,\ \zeta_2,\ \tau)-\nu_8(\zeta_1,\ \zeta_2,\ \tau)|$ \\
			\midrule
			
			0.5&	 0.1 & 4.237e-08 & 1.005e-11 & 1.554e-15 \\
			&	0.2 & 1.379e-06 & 1.302e-09 & 7.199e-13 \\
			&	0.3 & 1.065e-05 & 2.254e-08 & 2.796e-11 \\
			&	0.4 & 4.568e-05 & 1.710e-07 & 3.762e-10 \\
			&	0.5 & 1.419e-04 & 8.263e-07 & 2.832e-09 \\
			\midrule
			1&	0.1 & 8.474e-08 & 2.009e-11 & 3.109e-15 \\
			&	0.2 & 2.758e-06 & 2.605e-09 & 1.440e-12\\
			&	0.3 & 2.131e-05 & 4.508e-08 & 5.591e-11 \\
			&	0.4 & 9.136e-05 & 3.421e-07 & 7.524e-10 \\
			&	0.5 & 2.838e-04 & 1.653e-06 & 5.664e-09 \\
			
			\bottomrule
	\end{tabular}}
	\caption{Comparison of absolute errors for Example \ref{ex5}.}
	\label{tab:16}
\end{table}

\begin{figure}
	\centering
	\subfloat[$\nu_8(\zeta_1,\ \zeta_2, \ 0.5),\ \gamma=0.5$.]{%
		\resizebox*{6cm}{!}{\includegraphics{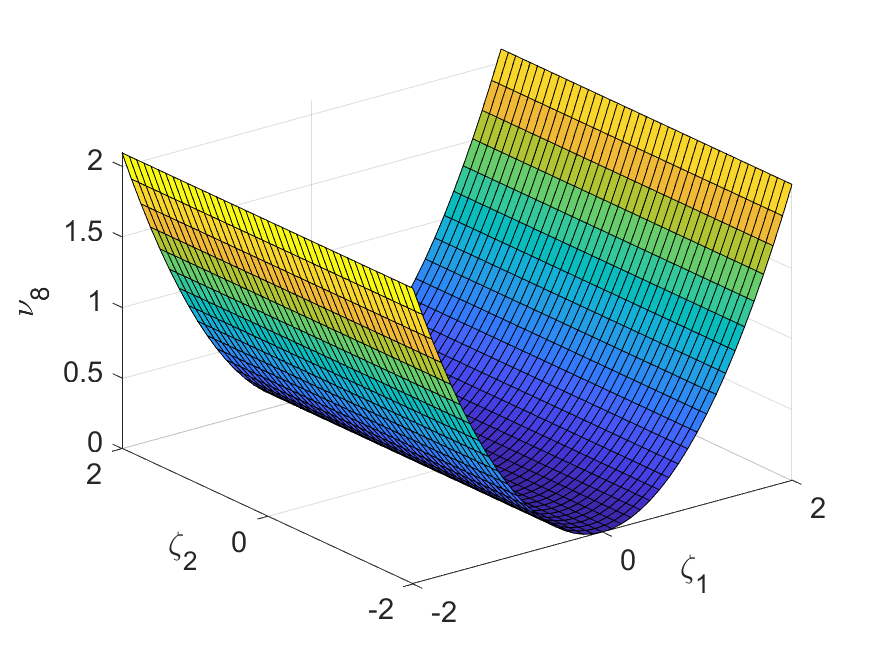}}}
	\subfloat[$\nu_8(\zeta_1,\ \zeta_2, \ 0.5),\ \gamma=0.75$.]{%
		\resizebox*{6cm}{!}{\includegraphics{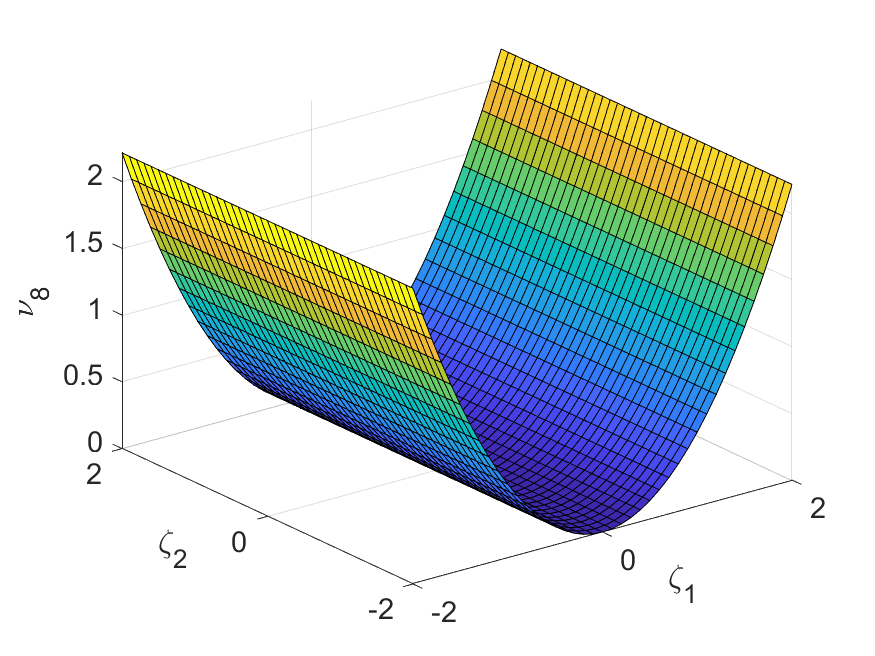}}}\hspace{5pt}
	\subfloat[$\nu_8(\zeta_1,\ \zeta_2, \ 0.5),\ \gamma=1$.]{%
		\resizebox*{6cm}{!}{\includegraphics{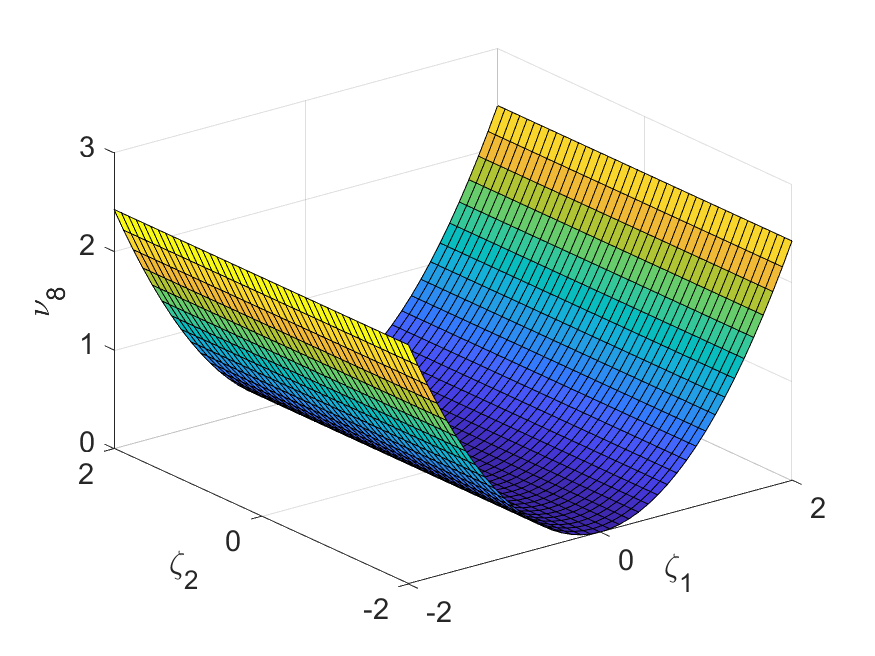}}}
	\subfloat[$\nu(\zeta_1,\ \zeta_2, \ 0.5),\ \gamma=1.$]{%
		\resizebox*{6cm}{!}{\includegraphics{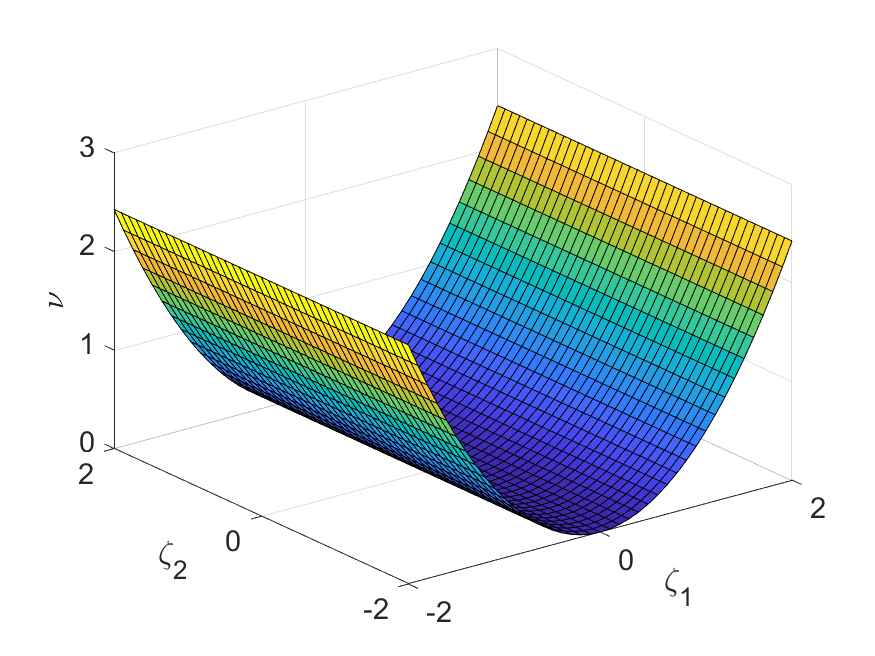}}}
	\caption{The variation of solution profile for various fractional orders for Example \ref{ex6}.} \label{fig_10}
\end{figure}

\begin{figure}
	\centering
	
	\subfloat[$|\nu(\zeta_1, \ \zeta_2, \ 0.5)-\nu_8(\zeta_1, \ \zeta_2, \ 0.5)|$ for $\gamma=0.5$.]{%
		\resizebox*{6cm}{!}{\includegraphics{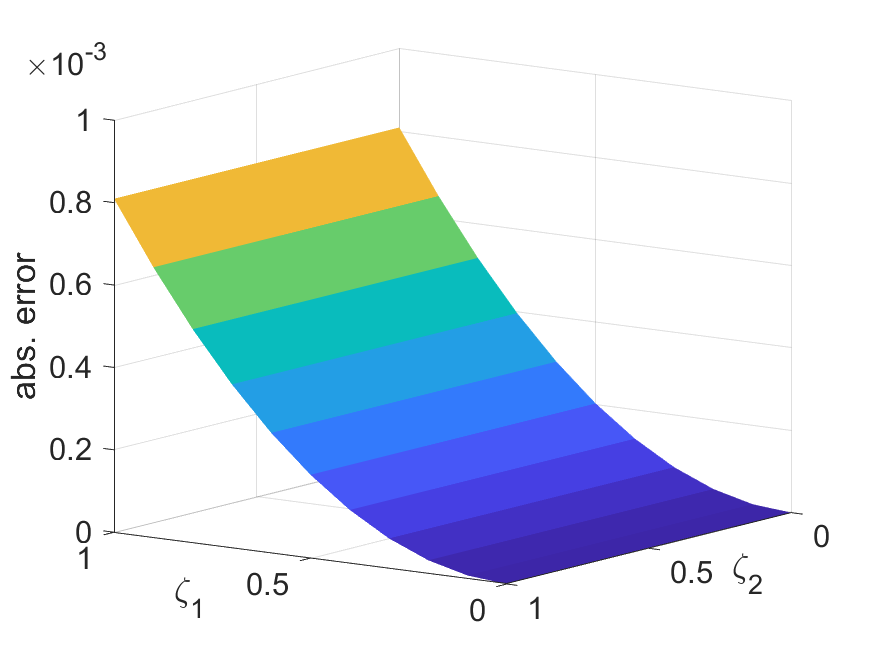}}} \hspace{1cm}
	\subfloat[$\zeta_1=\zeta_2=0.5,\ \gamma=1.$]{%
		\resizebox*{6cm}{!}{\includegraphics{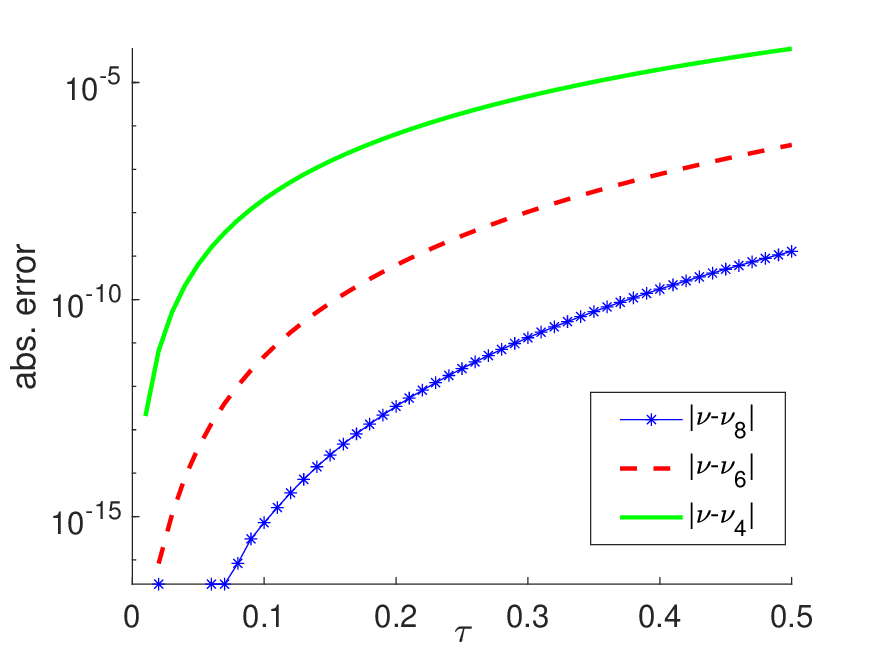}}}
	\caption{Error plots for Example \ref{ex6}.} \label{fig_11}
\end{figure}


\begin{table}
\centering
 \resizebox{16cm}{!}{
	\begin{tabular}{cccccc}
			\toprule
			$\tau$ & $\gamma=0.2$ & $\gamma=0.4$ & $\gamma=0.6$ & $\gamma=0.8$ & $\gamma=1$ \\
			\midrule
			0.15 & 0.144875205956596 & 0.159915115489348 & 0.179507802467414 & 0.198748036276498 & 0.215176994106291 \\
			0.30 & 0.143320696733757 & 0.142613098355304 & 0.154305440690060 & 0.168970871448092 & 0.185204555183594 \\
			0.45 & 0.148517386995271 & 0.132516456599877 & 0.137712605978767 & 0.146827012682291 & 0.159407038404212 \\
			0.60 & 0.157999864027977 & 0.126142617661293 & 0.125459902408442 & 0.129392823868808 & 0.137202915571429 \\
			0.75 & 0.170863152666341 & 0.122567529668948 & 0.115914327679565 & 0.115237877190844 & 0.118091686282839 \\
			0.90 & 0.186636215505243 & 0.121600198019816 & 0.108302496070487 & 0.103512796537055 & 0.101642659647879 \\
			\bottomrule
			
	\end{tabular}}
	\caption{LRPS solution $\nu_8(0.5,\ 0.5,\ \tau)$ for distinct $\gamma$ values for Example \ref{ex6}.}
	\label{tab:17}
\end{table}

\begin{table}
\centering
	{\begin{tabular}{ccccc}
			\toprule
			$\tau$ &exact & LRPS & abs. error & relative error \\
			\midrule
			0.15 & 0.215176994106264 & 0.215176994106291 & 2.609e-14 & 1.212e-13 \\
			0.30 & 0.185204555170429 & 0.185204555183594 & 1.316e-11 & 6.117e-11 \\
			0.45 & 0.159407037905443 & 0.159407038404212 & 4.987e-10 & 2.317e-09 \\
			0.60 & 0.137202909023507 & 0.137202915571429 & 6.547e-09 & 3.043e-08 \\
			0.75 & 0.118091638185254 & 0.118091686282839 & 4.809e-08 & 2.235e-07 \\
			0.90 & 0.101642414935150 & 0.101642659647879 & 2.447e-07 & 1.137e-06 \\
			\bottomrule

	\end{tabular}}
	\caption{Error estimates between $\nu(0.5, \ 0.5,\ \tau)$ and $\nu_8(0.5,\ 0.5,\ \tau)$ for $\gamma=1$ for Example \ref{ex6}.}
	\label{tab:18}
\end{table}

\begin{table}
\centering
	{\begin{tabular}{cccccc}
			\toprule
			$\zeta_1=\zeta_2$&	$\tau$ &$|\nu(\zeta_1,\ \zeta_2,\ \tau)-\nu_4(\zeta_1,\ \zeta_2,\ \tau)|$ & $|\nu(\zeta_1,\ \zeta_2,\ \tau)-\nu_6(\zeta_1,\ \zeta_2,\ \tau)|$ & $|\nu(\zeta_1,\ \zeta_2,\ \tau)-\nu_8(\zeta_1,\ \zeta_2,\ \tau)|$ \\
			\midrule
			0.5&	0.1 & 2.049e-08 & 4.899e-12& 7.216e-16 \\
			&		0.2 & 6.451e-07 & 6.194e-10 & 3.458e-13 \\
			&		0.3 & 4.820e-06 & 1.045e-08 & 1.316e-11 \\
			&		0.4 & 1.999e-05 & 7.738e-08 & 1.736e-10 \\
			&		0.5 & 6.004e-05 & 3.646e-07 & 1.281e-09 \\
			\midrule
			1&		0.1 & 8.196e-08 & 1.960e-11 & 2.887e-15 \\
			&		0.2 & 2.580e-06 & 2.478e-09 & 1.383e-12 \\
			&		0.3 & 1.928e-05 & 4.182e-08 & 5.266e-11\\
			&		0.4 & 7.995e-05 & 3.095e-07 & 6.945e-10\\
			&		0.5 & 2.402e-04& 1.458e-06 & 5.125e-09 \\
			
			\bottomrule
	\end{tabular}}
	\caption{Comparison of absolute errors for Example \ref{ex6}.}
	\label{tab:19}
\end{table}

\begin{figure}
	\centering
	\subfloat[$\nu_8(0.5, \ \zeta_2, \ \zeta_3, \ 0.5),\ \gamma=0.5$.]{%
		\resizebox*{6cm}{!}{\includegraphics{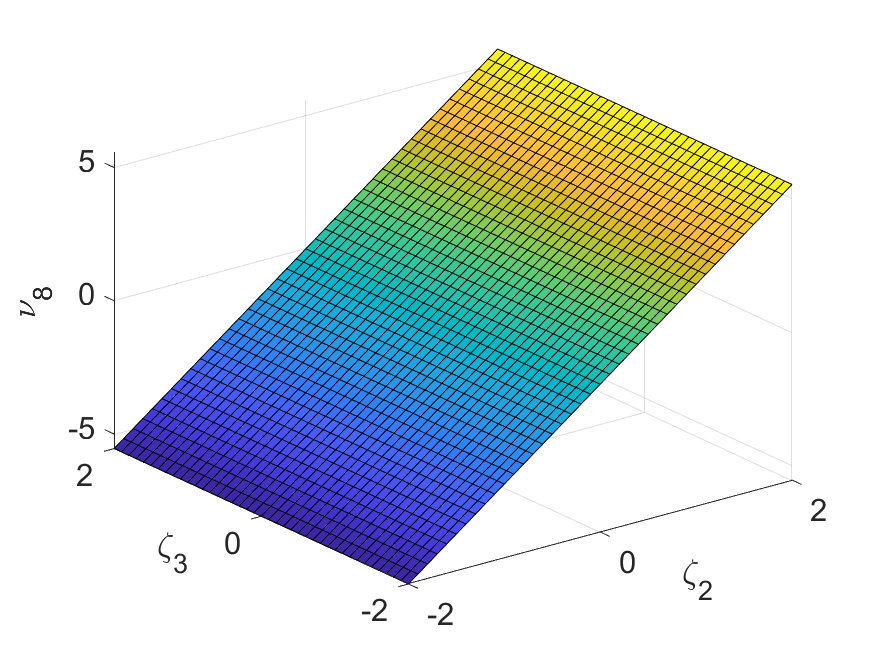}}}
	\subfloat[$\nu_8(0.5, \ \zeta_2, \ \zeta_3, \ 0.5),\ \gamma=0.75.$]{%
		\resizebox*{6cm}{!}{\includegraphics{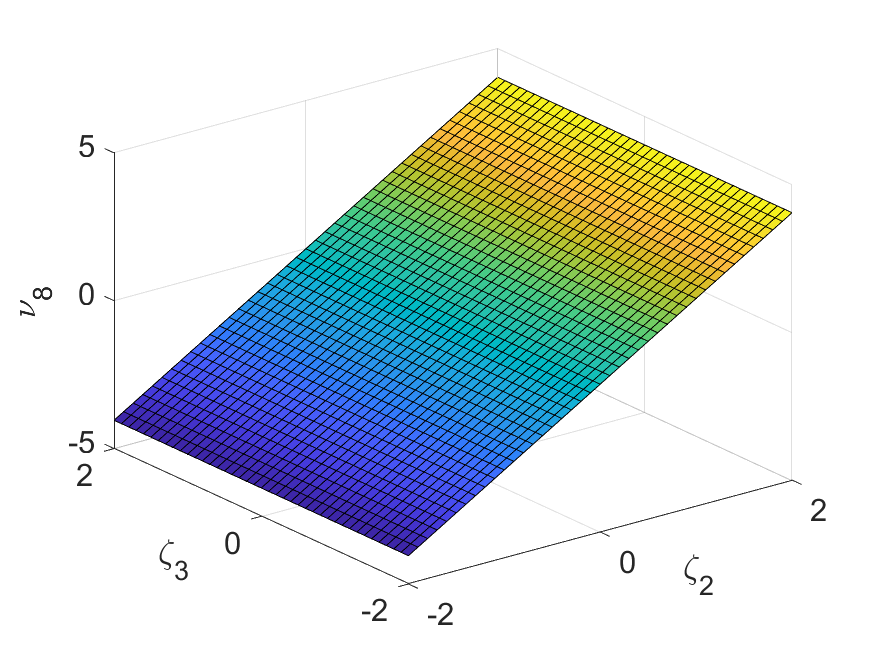}}}\hspace{5pt}
	\subfloat[$\nu_8(0.5, \ \zeta_2, \ \zeta_3, \ 0.5),\ \gamma=1$.]{%
		\resizebox*{6cm}{!}{\includegraphics{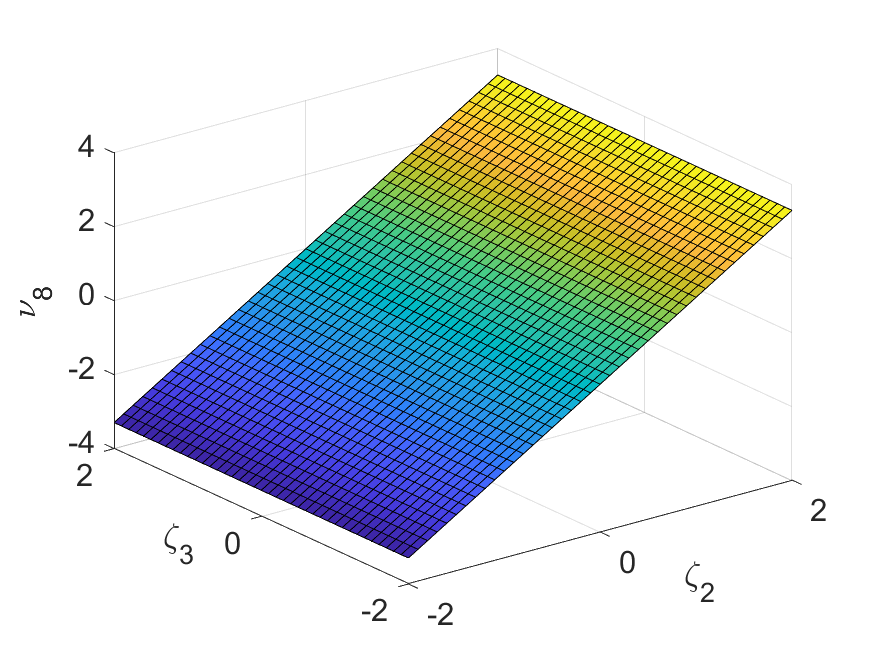}}}
	\subfloat[$\nu(0.5, \ \zeta_2, \ \zeta_3, \ 0.5),\ \gamma=1$.]{%
		\resizebox*{6cm}{!}{\includegraphics{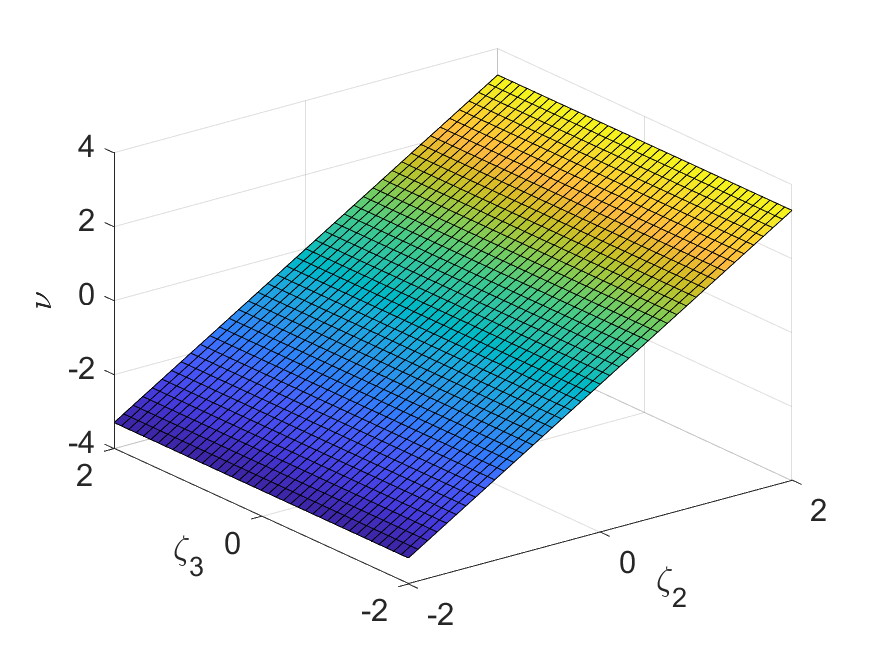}}}
	\caption{The variation of solution profile for various fractional orders for Example \ref{ex7}.} \label{fig_12}
\end{figure}

\begin{figure}
	\centering
	
	\subfloat[$|\nu(0.5,\ 0.5, \ \zeta_3,\tau)-\nu_8(0.5,\ 0.5,\  \zeta_3, \ \tau)|$ for $\gamma=0.5.$ 
    ]{%
		\resizebox*{6cm}{!}{\includegraphics{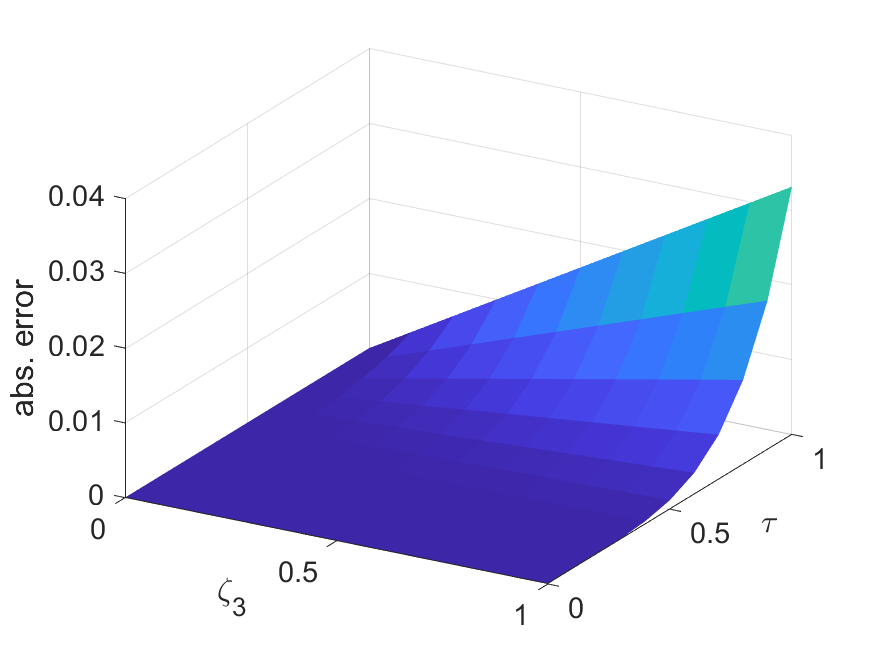}}} \hspace{1cm}
	\subfloat[$\zeta_1=\zeta_2=\zeta_3=0.5, \ \gamma=1$.]{%
		\resizebox*{6cm}{!}{\includegraphics{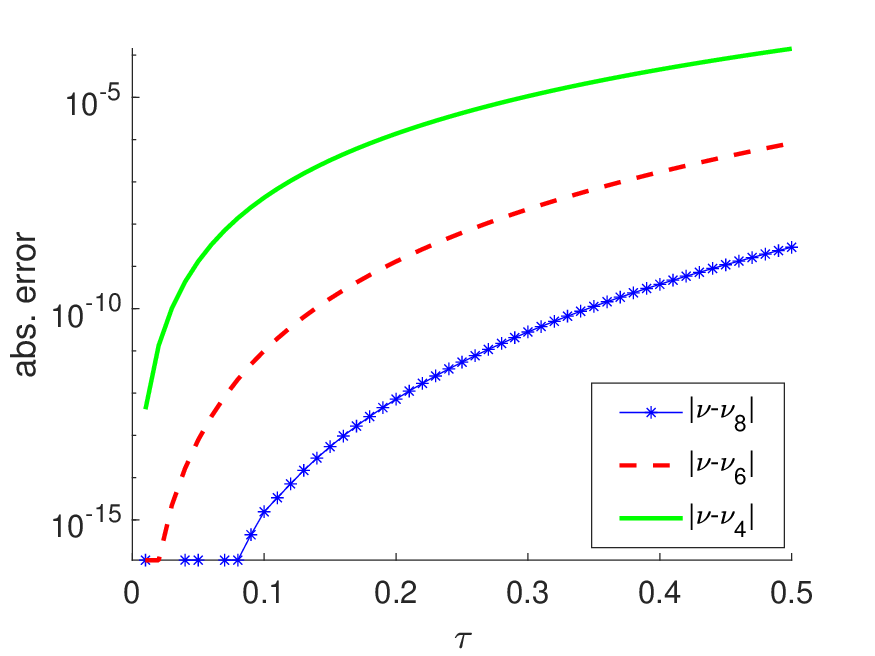}}}
	\caption{Error plots for Example \ref{ex7}.} \label{fig_13}
\end{figure}



\begin{table}
\centering
 \resizebox{16cm}{!}{
	\begin{tabular}{cccccc}
			\toprule
			$\tau$ & $\gamma=0.2$ & $\gamma=0.4$ & $\gamma=0.6$ & $\gamma=0.8$ & $\gamma=1$ \\
			\midrule
			0.15 & 1.600081043809555 & 0.952981800749650 & 0.737755085043236 & 0.636417247831176 & 0.580917121364088 \\
			0.30 & 2.135367406482500 & 1.247856716099542 & 0.925069682147884 & 0.766730117890084 & 0.674929403760045 \\
			0.45 & 2.631768336071708 & 1.556097169766161 & 1.126536612111816 & 0.911564290106330 & 0.784156091653569 \\
			0.60 & 3.116845586716595 & 1.893695854737707 & 1.352081318394315 & 1.076288836236974 & 0.911059385428571 \\
			0.75 & 3.599855147026588 & 2.268773048642961 & 1.608451623447892 & 1.265354801929119 & 1.058499896526337 \\
			0.90 & 4.084939830475697 & 2.687165687471149 & 1.901880218811209 & 1.483342806808060 & 1.229800969474331 \\
			\bottomrule
			
	\end{tabular}}
	\caption{LRPS solution $\nu_8(0.5,\ 0.5, \ 0.5,\ \tau)$ for distinct $\gamma$ values for Example \ref{ex7}.}
	\label{tab:20}
\end{table}


\begin{table}
\centering
	{\begin{tabular}{ccccc}
			\toprule
			$\tau$ &exact & LRPS & abs. error & relative error \\
			\midrule
			0.15 & 0.580917121364088 & 0.580917121364142 & 5.373e-14 & 1.000e-14 \\
			0.30 & 0.674929403760045 & 0.674929403788002 & 2.795e-11 & 4.224e-11 \\
			0.45 & 0.784156091653569 & 0.784156092745084 & 1.091e-09 & 1.127e-09 \\
			0.60 & 0.911059385428571 & 0.911059400195254 & 1.476e-08 & 1.634e-08 \\
			0.75 & 1.058499896526337 & 1.058500008306337 & 1.117e-07 & 1.056e-07 \\
			0.90 & 1.229800969474331 & 1.229801555578475 & 5.861e-07 & 5.747e-07 \\
			\bottomrule

	\end{tabular}}
	\caption{Error estimates between $\nu(0.5, \ 0.5,\ 0.5, \  \tau)$ and $\nu_8(0.5,\ 0.5,\ 0.5, \ \tau)$ for $\gamma=1$ for Example \ref{ex7}.}
	\label{tab:21}
\end{table}

\begin{table}
\centering
    \resizebox{16cm}{!}{
	{\begin{tabular}{cccccc}
			\toprule
			$\zeta_1=\zeta_2=\zeta_3$&	$\tau$ &$|\nu(\zeta_1, \ \zeta_2, \ \zeta_3, \ \tau)-\nu_4(\zeta_1, \ \zeta_2, \ \zeta_3, \ \tau)|$ & $|\nu(\zeta_1, \ \zeta_2, \ \zeta_3, \ \tau)-\nu_6(\zeta_1, \ \zeta_2, \ \zeta_3, \ \tau)|$ & $|\nu(\zeta_1, \ \zeta_2, \ \zeta_3, \ \tau)-\nu_8(\zeta_1, \ \zeta_2, \ \zeta_3, \ \tau)|$ \\
			\midrule
			
			0.5&	 0.1 & 4.237e-08 & 1.005e-11 & 1.554e-15 \\
			&	0.2 & 1.379e-06 & 1.302e-09 & 7.199e-13 \\
			&	0.3 & 1.065e-05 & 2.254e-08 & 2.796e-11 \\
			&	0.4 & 4.568e-05 & 1.710e-07 & 3.762e-10 \\
			&	0.5 & 1.419e-04 & 8.263e-07 & 2.832e-09 \\
			\midrule
			1&	0.1 & 8.474e-08 & 2.009e-11 & 3.109e-15 \\
			&	0.2 & 2.758e-06 & 2.605e-09 & 1.440e-12 \\
			&	0.3 & 2.131e-05 & 4.508e-08 & 5.591e-11 \\
			&	0.4 & 9.136e-05 & 3.421e-07 & 7.524e-10 \\
			&	0.5 & 2.838e-04 & 1.653e-06 & 5.664e-09 \\
			
			\bottomrule
	\end{tabular}}}
	\caption{Comparison of absolute errors for $\gamma=1$ for Example \ref{ex7}.}
	\label{tab:22}
\end{table}


\begin{figure}
	\centering
	\subfloat[$\nu_8(0.5,\ \zeta_2, \ \zeta_3, \ 0.5),\ \gamma=0.5$.]{%
		\resizebox*{6cm}{!}{\includegraphics{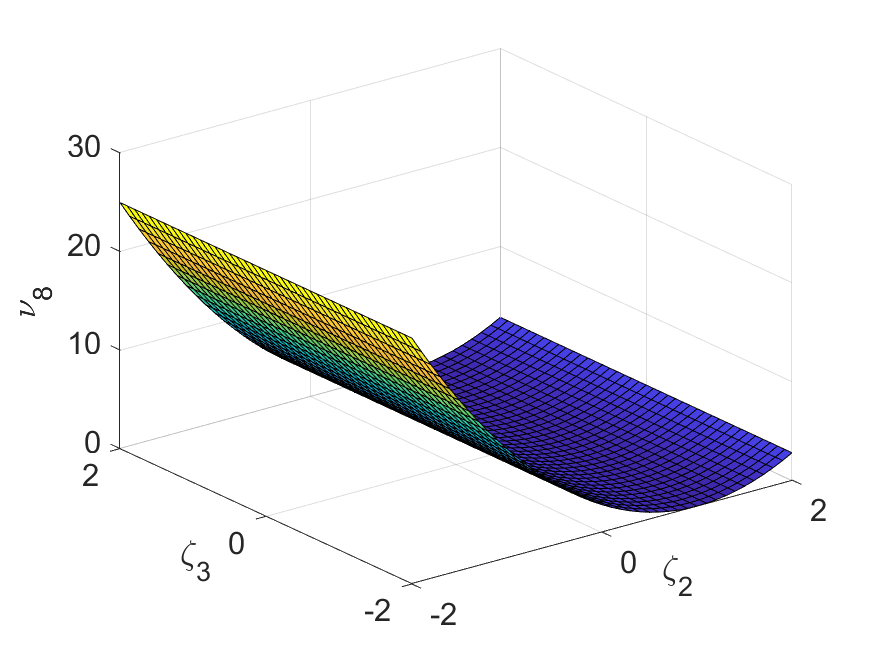}}}
	\subfloat[$\nu_8(0.5,\ \zeta_2, \ \zeta_3, \ 0.5),\ \gamma=0.75$.]{%
		\resizebox*{6cm}{!}{\includegraphics{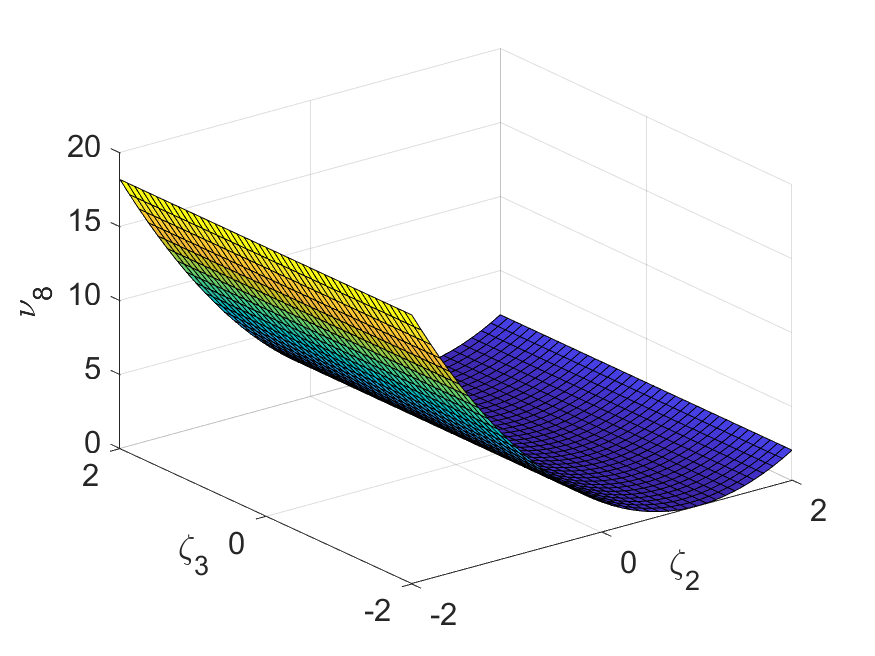}}}\hspace{5pt}
	\subfloat[$\nu_8(0.5,\ \zeta_2, \ \zeta_3, \ 0.5),\   \gamma=1$.]{%
		\resizebox*{6cm}{!}{\includegraphics{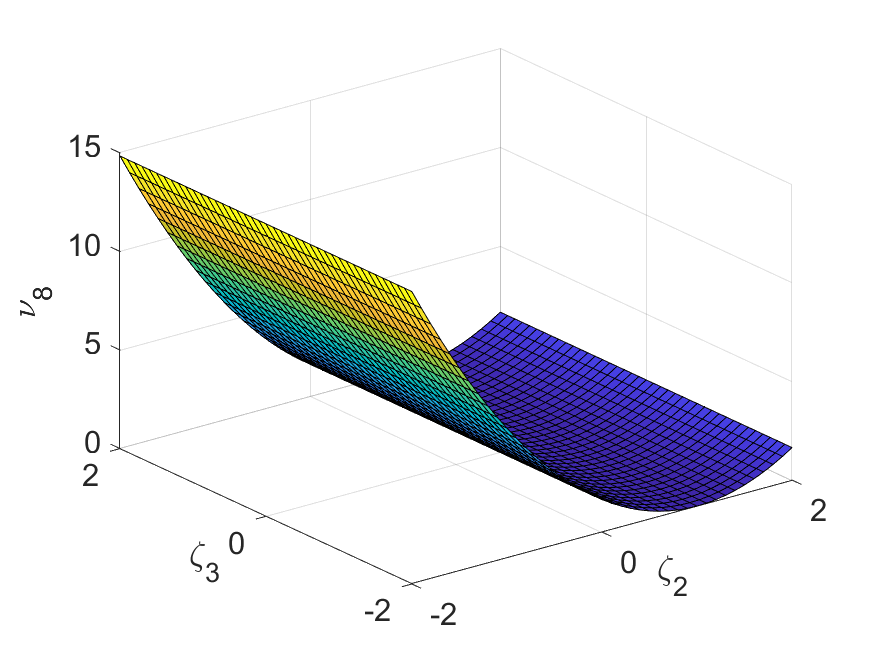}}}
	\subfloat[$\nu(0.5,\ \zeta_2, \ \zeta_3, \ 0.5),\  \gamma=1$.]{%
		\resizebox*{6cm}{!}{\includegraphics{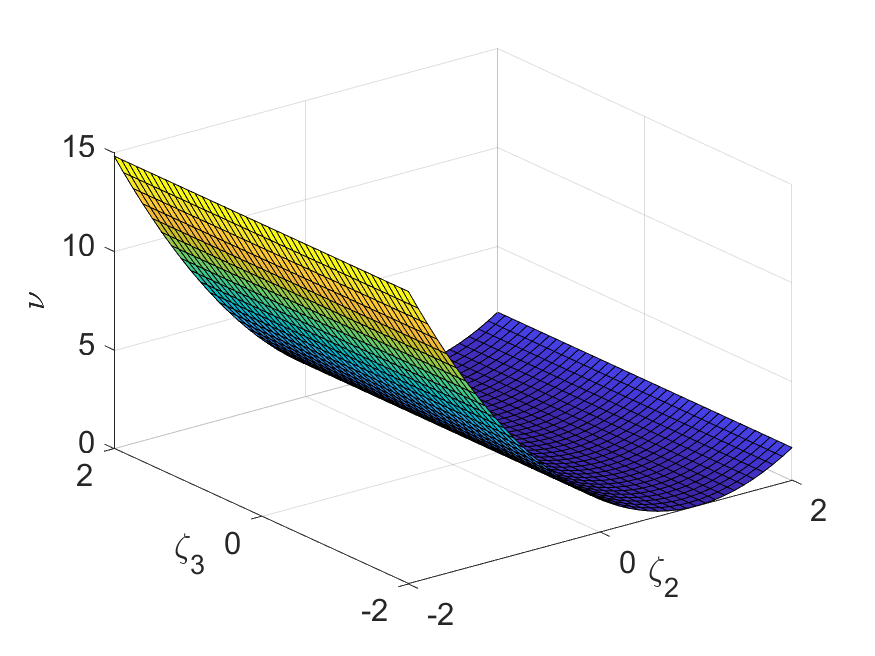}}}
	\caption{The variation of solution profile for various fractional orders for Example \ref{ex8}.} \label{fig_14}
\end{figure}

\begin{figure}
	\centering
	
	\subfloat[$|\nu(0.5,\ 0.5,\ \zeta_3, \ \tau)-\nu_8(0.5,\ 0.5,\ \zeta_3, \ \tau)|$ for $\gamma=0.5.$ 
    ]{%
		\resizebox*{6cm}{!}{\includegraphics{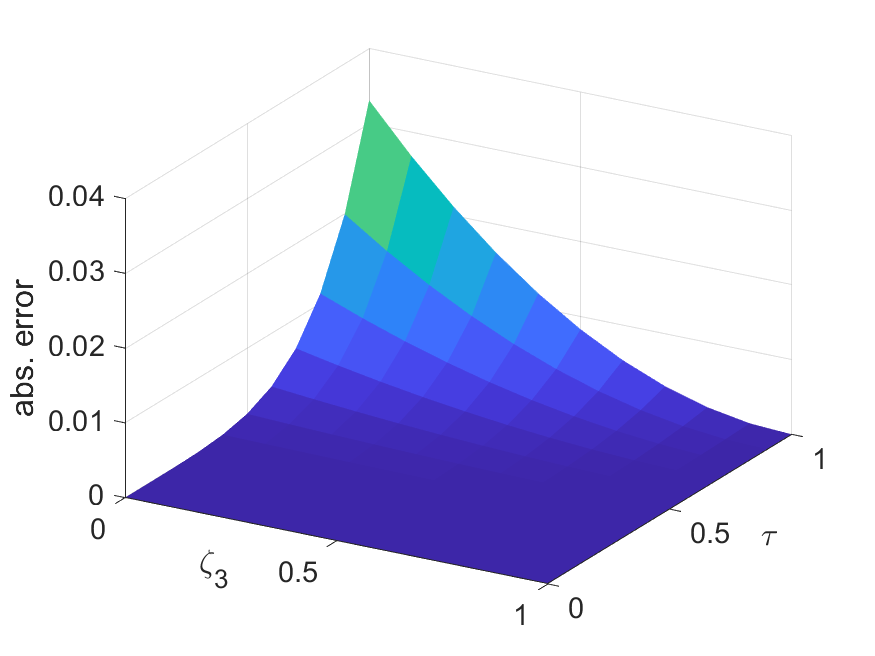}}} \hspace{1cm}
	\subfloat[$\zeta_1=\zeta_2=\zeta_3=0.5,\ \gamma=1$.]{%
		\resizebox*{6cm}{!}{\includegraphics{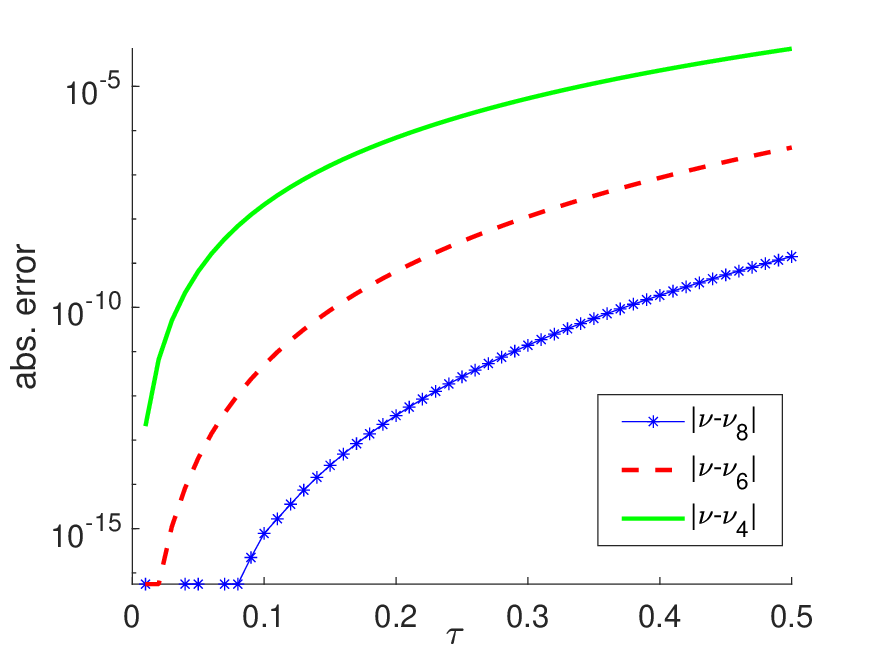}}}
	\caption{Error plots for Example \ref{ex8}.} \label{fig_15}
\end{figure}







\begin{table}
\centering
     \resizebox{16cm}{!}{
	\begin{tabular}{cccccc}
			\toprule
			$\tau$ & $\gamma=0.2$ & $\gamma=0.4$ & $\gamma=0.6$ & $\gamma=0.8$ & $\gamma=1$ \\
			\midrule
		0.15 & 0.800040521904777 & 0.476490900374825 & 0.368877542521618 & 0.318208623915588 & 0.290458560682044 \\
		0.30 & 1.067683703241250 & 0.623928358049771 & 0.462534841073942 & 0.383365058945042 & 0.337464701880022 \\
		0.45 & 1.315884168035854 & 0.778048584883080 & 0.563268306055908 & 0.455782145053165 & 0.392078045826784 \\
		0.60 & 1.558422793358298 & 0.946847927368853 & 0.676040659197157 & 0.538144418118487 & 0.455529692714286 \\
		0.75 & 1.799927573513294 & 1.134386524321481 & 0.804225811723946 & 0.632677400964559 & 0.529249948263168 \\
		0.90 & 2.042469915237849 & 1.343582843735574 & 0.950940109405604 & 0.741671403404030 & 0.614900484737165 \\
			\bottomrule
			
	\end{tabular}}
	\caption{LRPS solution $\nu_8(0.5,\ 0.5,\ 0.5,\ \tau)$ for distinct $\gamma$ values for Example \ref{ex8}.}
	\label{tab:23}
\end{table}


\begin{table}
\centering	
	{\begin{tabular}{ccccc}
			\toprule
			$\tau$ &exact  & LRPS & abs. error & relative error \\
			\midrule
			0.15 & 0.290458560682071 & 0.290458560682044 & 2.692e-14 & 9.269e-14 \\
			0.30 & 0.337464701894001 & 0.337464701880022 & 1.397e-11 & 4.142e-11 \\
			0.45 & 0.392078046372542 & 0.392078045826784 & 5.457e-10 & 1.391e-09 \\
			0.60 & 0.455529700097627 & 0.455529692714286 & 7.383e-09 & 1.621e-08 \\
			0.75 & 0.529250004153169 & 0.529249948263168 & 5.589e-08 & 1.056e-07 \\
			0.90 & 0.614900777789238 & 0.614900484737165 & 2.930e-07 & 4.765e-07 \\
			\bottomrule

	\end{tabular}}
	\caption{Error estimates between $\nu(0.5, \ 0.5,\ 0.5,\ \tau)$ and $\nu_8(0.5,\ 0.5,\ 0.5,\ \tau)$ for $\gamma=1$ for Example \ref{ex8}.}
	\label{tab:24}
\end{table}

\begin{table}
\centering
    \resizebox{16cm}{!}{
	{\begin{tabular}{cccccc}
			\toprule
			$\zeta_1=\zeta_2=\zeta_3$&	$\tau$ &$|\nu(\zeta_1, \ \zeta_2, \ \zeta_3, \ \tau)-\nu_4(\zeta_1, \ \zeta_2, \ \zeta_3, \ \tau)|$ & $|\nu(\zeta_1, \ \zeta_2, \ \zeta_3, \ \tau)-\nu_6(\zeta_1, \ \zeta_2, \ \zeta_3, \ \tau)|$ & $|\nu(\zeta_1, \ \zeta_2, \ \zeta_3, \ \tau)-\nu_8(\zeta_1, \ \zeta_2, \ \zeta_3, \ \tau)|$ \\
			\midrule
			
			0.5& 0.1   & 2.119e-08  & 5.023e-12  & 7.772e-16 \\
			&0.2   & 6.895e-07  & 6.512e-10  & 3.599e-13 \\
			&0.3   & 5.327e-06  & 1.127e-08  & 1.398e-11 \\
			&0.4   & 2.284e-05  & 8.552e-08  & 1.881e-10 \\
			&0.5   & 7.094e-05  & 4.132e-07  & 1.416e-09 \\
			\midrule
			0.75	&0.1   & 5.296e-09  & 1.256e-12  & 1.943e-16 \\
			&0.2   & 1.724e-07  & 1.628e-10  & 8.998e-14 \\
			&0.3   & 1.332e-06  & 2.817e-09  & 3.495e-12 \\
			&0.4   & 5.710e-06  & 2.138e-08  & 4.702e-11 \\
			&0.5   & 1.774e-05  & 1.033e-07  & 3.540e-10 \\
			\bottomrule
	\end{tabular}}}
	\caption{Comparison of absolute errors 
for Example \ref{ex8}.}
	\label{tab:25}
\end{table}


This section analyzes the numerical results and graphical illustrations of numerical problems 
solved in the previous section.
The surface plots depict 
the variation of approximate solutions according to various fractional orders. The high accuracy of this method is highlighted through a comparison analysis of the exact and the LRPS solution. A comparative analysis is performed to validate our results, comparing the current method with the existing method, and the findings are tabulated.


The approximate analytical solution $\nu_8(\vec{\zeta},\tau)$ of TFFPE obtained from the LRPS method is depicted through 3-dimensional plots.  Figures \ref{fig_1}, \ref{fig_2}, \ref{fig_6}, \ref{fig_8}, \ref{fig_10}, \ref{fig_12}, and \ref{fig_14} represent the  3-dimensional plots of the LRPS solution $\nu_8(\vec{\zeta},\tau)$ and the exact solution $\nu(\vec{\zeta},\tau)$ for various fractional orders 
for Examples \ref{ex1} - \ref{ex8} respectively. 
It is notable that the approximate solution $\nu_8(\vec{\zeta},\tau)$ aligns with the exact solution as $\gamma$ tends to 1. There is a good harmony between the approximate solution $\nu_8(\vec{\zeta},\tau)$  and the exact solution $\nu(\vec{\zeta},\tau)$.
Figures \ref{fig_1a}, \ref{fig_2a}, and \ref{fig_4a} illustrate the fluctuation of the approximate solution $\nu_8(\zeta,\tau)$ at different times $\tau$ with the constant fractional order $\gamma=0.5.$  Moreover, it is clear that as time $\tau$ progresses, the solution grows gradually. 

As the exact solution of all considered numerical problems of TFFPE for $\gamma=1$ is known, we utilize various error measures for validating the efficacy of the employed method, 
which are defined as follows: 
\vskip 0.225cm
\begin{itemize}
	\item []
	Abs. error $= |\nu(\vec{\zeta},\tau)-\nu_n(\vec{\zeta},\tau)|$
	\item []	Relative error $=\bigg|\frac{\nu(\vec{\zeta},\tau)-\nu_n(\vec{\zeta},\tau)}{\nu(\vec{\zeta},\tau)}\bigg|$
\end{itemize} 
The absolute error plots are represented in Figures \ref{fig_3},  \ref{fig_7}, \ref{fig_9}, \ref{fig_11}, \ref{fig_13}, and \ref{fig_15} for Examples \ref{ex1} - \ref{ex8} respectively. The subplots (a) are the 3-dimensional plots of absolute error (abs. error), $|\nu(\vec{\zeta},\tau)-\nu_8(\vec{\zeta},\tau)|$ for $\gamma=0.5.$ The subplots (b) are the 2-dimensional plots of different absolute errors while considering approximate solutions $\nu_4,$ $\nu_6,$ and $\nu_8$ on the logarithmic scale.  The absolute error is closely aligned with zero if we consider $\nu_8(\vec{\zeta},\tau).$ We observe that the absolute error grows with an increase in $\tau$.  Furthermore, it is evident that an increased number of terms in the series solution correlates with enhanced accuracy and reduced absolute error.

The variation of solution value $\nu_8(\vec{\zeta},\tau)$ for various fractional orders is tabulated in Tables \ref{tab:1}, \ref{tab:3}, \ref{tab:11}, \ref{tab:14}, \ref{tab:17}, \ref{tab:20}, and \ref{tab:23}, for Examples \ref{ex1} - \ref{ex8}, respectively. 
 These tables show the influence of fractional orders on the behavior of the solution. As the fractional order tends to 1, the solution increases with time for all examples except example \ref{ex6}. But in Example \ref{ex6}, 
 , the solution exhibits an increasing trend with increasing $\gamma$ for small values of 
 $\tau$. However, the solution shows a decreasing trend as time progresses particularly for larger $\tau$ values, especially evident when $\gamma$ is closer to 1.

 The absolute and relative errors for the case of $\gamma=1$ for Examples \ref{ex2}-\ref{ex8} are summarized in Tables \ref{tab:5},  \ref{tab:12}, \ref{tab:15}, \ref{tab:18}, \ref{tab:21}, and \ref{tab:24}, respectively. These tables indicate that the absolute errors decrease as time decreases, while they increase as time progresses. We observe that the LRPS solution is consistent with the exact solution by taking only eight terms. Conclusively, the absolute errors and the relative errors are within an admissible range.
 Tables \ref{tab:6}, \ref{tab:13}, \ref{tab:16}, \ref{tab:19}, \ref{tab:22}, and \ref{tab:25} are the comparative study of absolute errors for $k=4,$ $k=6,$ and $k=8$ for the Examples \ref{ex2} - \ref{ex8}, respectively for $\gamma=1$. These tables clearly illustrate that, as we incorporate more terms in the approximate solution, we achieve improved results with reduced errors.
 To validate our outcomes,  a comparison of the suggested method with other existing methods has been conducted and tabulated in Tables \ref{tab:2} and \ref{tab:4}. The derived solution is consistent with the existing method \cite{HAS24} for Example \ref{ex1}. 
  Thus, the comparative analysis with existing methods confirms that the proposed approach is effectively utilized to solve TFFPE.

 In this manuscript, the LRPS method is successfully applied to tackle both linear and non-linear TFFPE in multidimensional space. 
 The following are the key implications and limitations  
 of the proposed method.
\begin{itemize}
	\item This method offers a highly accurate semi-analytical solution for fractional differential equations, providing a deeper analysis of the system than numerical methods.
	\item This is a versatile method that is applicable to an extensive spectrum of linear or non-linear fractional differential equations. The derivative in the spatial direction can be easily tackled using the Laplace transform.
	\item Higher dimensional fractional differential equations can be easily addressed and solved by utilizing this approach.
	\item This method enables us to enhance the precision of the solution by adding the number of terms in the series. More number of terms produce improved outcomes.
	\item The method is simple to execute, necessitating only fundamental operations and straightforward calculations, which could benefit computational efficiency.
	\item This technique does not require discretization, perturbation, or linearization. Moreover, this method is superior to the classical residual power series method, as there is no need to compute the fractional derivatives.
	\item  This method typically needs a solution that has the power series form.  A lot of fractional differential equations lack closed-form solutions. The reliance on series expansions may not always yield useful approximations, especially if the series converges poorly.

\end{itemize}  
\section{Consequence of the control term in TFFPE }
The LRPS technique is an effective and precise approach that we successfully used to tackle the non-linear, multi-dimensional TFFPEs in the absence of a control term.  Now, the question arises about the consequence of control terms in the numerical simulation of TFFPEs. We consider the one-dimensional TFFPE with a control term
\begin{equation}
    D^\tau _\gamma \nu(\zeta, \tau) = D_\zeta \left( \frac{\zeta}{2} \nu(\zeta, \tau)\right) +D_{\zeta \zeta} \nu (\zeta, \tau) + E_\gamma (\tau^ \gamma), \hspace{0.2cm} (\zeta) \in\mathbb{R},\hspace{0.2cm} \tau \geq 0.
\end{equation} subject to the initial condition $\nu (\zeta,0)= \zeta+2$. The exact solution for this problem is $\nu(\zeta, \tau) = (\zeta +2) E_\gamma (\tau^\gamma)$. After applying the discussed strategy to this TFFPE, we obtain the solution in the infinite series form as: 
\begin{equation}
    \nu(\zeta, \tau)= \sum_{k=0} ^\infty \frac{(\zeta +2) \tau ^{k\gamma}}{\Gamma (k\gamma +1)}.
\end{equation} Thus, our approach provides the solution which converges to the exact solution for all $\gamma \in (0,1]$. 

 However, this method may not be suitable for all TFFPEs that have a control term. Assume the control term includes functions such as $\frac{1}{\tau}, ~ \exp{\tau ^2}$, or other functions for which the LT 
 cannot be defined. Consequently, we fail to use the LRPS approach in these circumstances. Moreover, even if the control terms are functions with a specified LT, 
 the LRPS method may not provide a solution.

For example, consider the one-dimensional TFFPE, 
\begin{equation}
    D_\tau ^\gamma \nu (\zeta, \tau) = D_\zeta \left(\exp{(\zeta-0.5)^2} \nu(\zeta, \tau)\right)+D_{\zeta \zeta} \nu (\zeta,\tau) + \mathfrak{g}(\zeta, \tau), ~~ 0<\zeta<1,~ 0<\tau \leq 1, 
\end{equation} with the initial condition $\nu(\zeta, 0) = 0$ and the control term  \begin{equation*}
    \mathfrak{g}(\zeta, \tau) = \frac{\Gamma (3)}{\Gamma(3-\gamma)} \tau ^{2-\gamma} \sin{(\pi \zeta)} - \tau ^2 \exp{(\zeta-0.5)^2} \left ( 
2(\zeta-0.5)\sin{(\pi \zeta)} +\pi \cos{(\pi \zeta)} \right) + \pi^2\tau ^2 \sin{(\pi \zeta)}.
\end{equation*} The exact solution for this TFFPE is $\nu(\zeta, \tau) =\tau ^2 \sin{(\pi \zeta)}$. 

 After applying LT and utilising the initial condition, we obtain the $ k^ {th}$ LRF as, 
 \begin{equation} \label{C68}
 \begin{split}
     \mathcal{LR}_k(\zeta, s) = &~ \mathcal{V}_k(\zeta,s)  -\frac{1}{s^\gamma} D_\zeta \left ( \exp{(\zeta-0.5)^2} \mathcal{V}_k(\zeta, s) \right) - \frac{1}{s^\gamma}D_{\zeta \zeta} \mathcal{V}_k(\zeta,s) - \frac{2\sin{(\pi \zeta)}}{s^3} \\
     &- \frac{2}{s^{3+\gamma}} \left(  \pi ^2 \sin{(\pi \zeta)} - \exp{(\zeta -0.5)^2}\left[ 2(\zeta-0.5)\sin{(\pi \zeta)+\pi \cos{(\pi \zeta)}} \right] \right)
 \end{split}
 \end{equation}
We replace the series representation of $\mathcal{V}_k(\zeta,s)$ in Eq. \eqref{C68} and, as previously done, impose the limit condition. For \( k=1 \), the coefficient is obtained as \( p_1(\zeta) = 0 \). For \( k \geq 2 \), we have, 
\begin{equation}\begin{split}   
     s^{2\gamma +1}\mathcal{LR}_2(\zeta, s) &= p_2 (\zeta) - \frac{ p_2^{``}(\zeta)}{s^{\gamma}}- \frac{1}{s^\gamma}\left(  2(\zeta-0.5)\exp{(\zeta-0.5)^2}p_2(\zeta) + \exp{(\zeta-0.5)^2}p_2'(\zeta) \right) - \frac{2\sin{(\pi \zeta)}}{s^{2-2\gamma}} \\ & -\frac{1}{s^{2-\gamma}}   \left(  \pi ^2 \sin{(\pi \zeta)} - \exp{(\zeta -0.5)^2}\left[ 2(\zeta-0.5)\sin{(\pi \zeta)+\pi \cos{(\pi \zeta)}} \right] \right).     
\end{split}
\end{equation}
Thus, we examine two distinct cases of \( \gamma \): \( 0 < \gamma < 1 \) and \( \gamma = 1 \).  For $0< \gamma <1$, the limit condition $\lim_{s\rightarrow \infty} s^{2\gamma + 1}\mathcal{LR}_2(\zeta, s) =0$ gives the coefficient $p_2(\zeta) =0.$ 
  However, for $\gamma = 1$, we get $p_2(\zeta) = 2\sin{(\pi \zeta)}$. Now, we designate \( k = 3 \). Thus, for $0<\gamma <1,$
  \begin{equation}
      \begin{split}
          s^{3\gamma +1}\mathcal{LR}_3(\zeta, s) &=p_3(\zeta) -\frac{p_3''(\zeta)}{s^\gamma}-\frac{1}{s^\gamma}\left(  2(\zeta-0.5)\exp{(\zeta -0.5)^2}p_2(\zeta) + \exp{(\zeta-0.5)^2}p_2'(\zeta) \right)
          -\frac{2\sin{(\pi \zeta)}}{s^{2-3\gamma}}\\ 
          &- -\frac{1}{s^{2-2\gamma}}   \left(  \pi ^2 \sin{(\pi \zeta)} - \exp{(\zeta -0.5)^2}\left[ 2(\zeta-0.5)\sin{(\pi \zeta)+\pi \cos{(\pi \zeta)}} \right] \right). 
      \end{split}
  \end{equation}Hence, the expression \( s^{3\gamma + 1} \mathcal{LR}_3(\zeta,s) \) diverges to infinity as \( s \rightarrow \infty \), specifically for \( \gamma \in \left(\frac{2}{3}, 1\right) \). Thus, this violates condition $\lim_{s \rightarrow \infty} s^{3\gamma +1}\mathcal{LR}_3(\zeta,s)=0$, making it impossible to estimate the coefficient $p_3(\zeta)$ for $\frac{2}{3} < \gamma <1$. For $\gamma \notin (\frac{2}{3},1)$, we find that $p_3(\zeta) = 0$, but for $\gamma = \frac{2}{3}$, $p_3(\zeta) = 2\sin{(\pi \zeta)}$. Consequently, the suggested approach does not provide an accurate solution consistently in the present situation for all $\gamma \in (0,1]$.

Therefore, the implementation of the proposed technique is dependent on the control term included in the TFFPE. The proposed approach may provide an accurate and efficient solution for TFFPE with a control term.
\section{Conclusion}
The time-fractional Fokker-Planck equation is an extension of the classical Fokker-Planck equation, which has great physical significance in the real world. 
This study develops a semi-analytical solution for the time-fractional Fokker-Planck problem by employing the Laplace residual power series approach, a blend of the Laplace transform with the residual power series method. The potency of the Laplace transform converts the target fractional differential equation into a simple algebraic equation that is simpler to solve than solving a fractional differential equation. We solve a range of linear and non-linear numerical problems in multi-dimensional space to show the validity and capability of the method.  Furthermore, the analysis features a discussion on error evaluation, and the characteristics of the resulting solution offers insights into the method's reliability. The influence of fractional order on the solution is displayed graphically. The outcomes are in excellent agreement with the existing techniques. The implementation of the suggested approach for time-fractional Fokker-Planck equations with a control term is contingent upon the control term involved.  
   
   Thus, the Laplace residual power series approach is a reliable and effective tool for solving a wide range of fractional differential equations that characterize chaotic dynamics, anomalous transport, solitonic wave propagation, and nonlinear phenomena. Future studies could include the extension of the application of the proposed method for the space-time fractional Fokker-Planck equation. 

\section*{Acknowledgments}
The authors would like to extend their gratitude to anonymous editors and reviewers for their suggestions.

\section*{Financial disclosure}
First author is supported with faculty research grant by NIT Calicut. The second author is supported by University Grant Commission of India (Ref No. 221610030044).
\section*{Data availability}
Available from the corresponding author upon reasonable request.

\section*{Author contributions statement}

N.G.: PhD supervisor, review and edit the manuscript. V.R.: computational study, write the main manuscript. 
\section*{Declarations}
\section*{Competing interests}
The authors declare no competing interests.

\section*{Authors and Affiliations}
Dr. Neetu Garg, Faculty, Department of Mathematics, National Institute of Technology, Calicut, 673601 Kerala, India\\
Varsha R., Research Scholar, Department of Mathematics, National Institute of Technology, Calicut, 673601 Kerala, India

\section*{Additional information}
Correspondence should be addressed to Neetu Garg.
\end{document}